\documentclass[11pt,twoside]{amsart}

\usepackage{geometry}
\usepackage{amsmath,amssymb,amsthm}
\title[Fully nonlinear CR invariant equations on $\h^n$]{On fully nonlinear CR invariant equations on the Heisenberg group}

\email{yyli@math.rutgers.edu (Yanyan Li)}
\email{dario.monticelli@gmail.com (Dario Daniele Monticelli)}



\newcommand{\e}{\varepsilon}
\newcommand{\Hn}[1]{|#1|_{\!_H}}
\newcommand{\Hgrad}[1]{\nabla_{\!_H}#1}
\newcommand{\Hhesss}[1]{\nabla_{\!_{H,s}}^2#1}
\newcommand{\field}[1]{\mathbf{#1}}
\newcommand{\R}{\field{R}}
\newcommand{\h}{\field{H}}
\newcommand{\C}{\field{C}}
\newcommand{\N}{\field{N}}

\newcommand{\crochet}[2]{\left\langle{#1},\,{#2}\right\rangle}

\newcommand{\wt}{\widetilde}

\newtheorem{thm}{\textbf{Theorem}}[section]
\newtheorem{lem}[thm]{\textbf{Lemma}}
\newtheorem{pro}[thm]{\textbf{Proposition}}
\newtheorem{cor}[thm]{\textbf{Corollary}}
\newtheorem{defn}[thm]{\textbf{Definition}}
\theoremstyle{remark}
\newtheorem{rem}[thm]{\textbf{Remark}}

\newtheoremstyle{Claim}{}{}{\itshape}{}{\itshape\bfseries}{:}{ }{#1}
\theoremstyle{Claim}

\begin{document}

\maketitle \vspace{-0,3cm}

\begin{center}
\textsc{\textmd{Y. Y. Li\footnote{Rutgers University, New Brunswick, NJ, USA} and
D. D. Monticelli\footnote{Universit\`{a} degli Studi di Milano, Italy. Partially
supported by GNAMPA, project ``Equazioni differenziali e sistemi
dinamici''. Partially supported by MIUR, project ``Metodi
variazionali e topologici nello studio di fenomeni nonlineari''.}}}
\end{center}


\begin{abstract}
In this paper we provide a characterization of second order fully
nonlinear CR invariant equations on the Heisenberg group, which is
the analogue in the CR setting of the result proved in the Euclidean
setting by A. Li and the first author in \cite{LiLi}. We also prove
a comparison principle for solutions of second order fully nonlinear
CR invariant equations defined on bounded domains of the Heisenberg
group and a comparison principle for solutions of a family of second
order fully nonlinear equations on a punctured ball.
\end{abstract}

\section{\textbf{Introduction and main results}}\label{intro}

As it was pointed out by Jerison and Lee in \cite{JerLee1}, there
are important similarities between conformal geometry and the
geometry of CR manifolds, which serve as abstract models of real
hypersurfaces in complex manifolds, which we are going to discuss
briefly in the following. On this subject see also the survey
article by Beals, Fefferman and Grossman \cite{BeaFefGro} and the
book by Dragomir and Tomassini \cite{TomDra}.

A CR manifold is a differentiable manifold $M$ equipped with a
subbundle $\mathcal{H}$ of the complexified tangent bundle $\C TM =
TM \otimes \C$ such that $[\mathcal{H},\mathcal{H}]\subseteq
\mathcal{H}$ (i.e. $\mathcal{H}$ is formally integrable) and
$\mathcal{H}\cap\bar{\mathcal{H}}=\{0\}$ (i.e. $\mathcal{H}$ is
almost Lagrangian). The bundle $\mathcal{H}$ is called a CR
structure on the manifold $M$. An abstract CR manifold is said to be
of hypersurface type if $\text{dim}_{\R}M=2n+1$ and
$\text{dim}_{\C}\mathcal{H}=n$.


If the CR manifold M is of hypersurface type and oriented, it is
possible to associate to its CR structure $\mathcal{H}$ a one form
$\theta$ globally defined on $M$ such that
$Ker(\theta)=\mathcal{H}$. $\theta$ is unique modulo a multiple of
nonzero function on $M$:  a choice of a nonzero multiple of such
$\theta$ is called a pseudohermitian structure on $M$.   $d\theta$
defines the Levi form $L_\theta$ on $\mathcal{H}$. If the Levi form
$L_\theta$ is strictly positive definite, we say that $M$, with this
$CR$ structure, is strictly pseudoconvex. In this case, the form
$\theta$, as well as a nonzero function multiple of $\theta$, is a
contact form on $M$.  See for instance \cite{TomDra}.

A scalar curvature associated to the pseudohermitian structure has
been introduced by Webster in \cite{Web1}, \cite{Web2} and by Tanaka
in \cite{Tan}. Thus the CR Yamabe problem is formulated as follows:
on a strictly pseudoconvex CR manifold, find a choice of
pseudohermitian structure, or equivalently a choice of contact form
$\theta$, with constant pseudohermitian scalar curvature.

Jerison and Lee proved in \cite{JerLee1} that there exists a CR
numerical invariant $\lambda(M)$ associated to every compact
strictly pseudoconvex orientable CR manifold $M$ of dimension
$2n+1$, such that $\lambda(M)$ is always less than or equal to the
value corresponding to the sphere $S^{2n+1}\subset\C^{n+1}$ and such
that the CR Yamabe problem admits solution on $M$, provided that
$\lambda(M)<\lambda(S^{2n+1})$. This result is an analogue in the CR
setting of the classical result of Aubin \cite{Aub} on Riemannian
manifolds.

Jerison and Lee also proved in \cite{JerLee2} that if $\theta$ is a
contact form associated to the standard CR structure on $S^{2n+1}$
having constant pseudohermitian scalar curvature, then it is
obtained from a constant multiple of the standard contact form on
$S^{2n+1}$ via a CR automorphism of the sphere. This result is then
an analogue in the CR setting of the well known result by Obata in
\cite{Oba} and by Gidas, Ni and Nirenberg in \cite{GidNiNi}.

The Heisenberg group $\h^n$ is CR equivalent to the sphere
$S^{2n+1}\subset\C^{n+1}$ minus a point via the Cayley transform,
see e.g. \cite{JerLee1}, so that the Heisenberg group plays in CR
geometry the same role as $\R^n$ in conformal geometry while the
Cayley transform corresponds to the stereographic projection.

Under a change of pseudohermitian structure given by
$\wt{\theta}=u^{p-2}\theta$ with $p=\frac{2n+2}{n}$, the
pseudohermitian scalar curvature $R$ changes according to the
following equation
\begin{equation}\label{58}
b_n\Delta_bu + Ru = \wt{R}u^{p-1},\qquad b_n:=\frac{2n+2}{n},
\end{equation}
as one can find in \cite{JerLee1}. Here $\Delta_b$ is the
sublaplacian operator on the CR manifold $M$, which is a linear
second order subelliptic operator, see also \cite{TomDra},
\cite{JerLee1} and references therein. On the sublaplacian on the
Heisenberg group, which we will denote by $\Delta_{\!H}$, see also
section \ref{notation}.

One interesting feature of equation \eqref{58} is that the exponent
$p$ in the nonlinearity is the same as the one appearing in a
Sobolev--type inequality for functions in
$\mathcal{C}_0^\infty(\h^n)$, which is related to the CR structure
defined on $\h^n$, that was proved by Folland and Stein in
\cite{FolSte2}.

We report here a very nice table that summarizes many important
similarities between CR geometry and conformal geometry, as it
appears in \cite{JerLee1}.\\

\begin{tabular}{ll}
  \vspace{0,3cm}\underline{Conformal geometry} & \underline{CR geometry} \\
  Riemannian Manifold $(M,g)$\hspace{0.8cm} & CR manifold $(M,\theta)$ \\
  Euclidean space $\R^m$ & Heisenberg group $\h^n$ \\
  $m$--sphere $S^{m}\subset\R^{m+1}$  & $(2n+1)$--sphere $S^{2n+1}\subset\C^{n+1}$\\
  Stereographic projection & Cayley transform \\
  Riemannian normal coordinates & Folland-Stein normal coordinates\\
  Scalar curvature K & Webster (pseudohermitian) scalar curvature R\\
  Laplace-Beltrami operator $\Delta$ & Sublaplacian operator $\Delta_b$\\
  Sobolev spaces $W^{k,r}$ & Folland-Stein spaces $S^{k,r}$\\
  Sobolev embedding $W^{1,2}\subset L^q$,  & Sobolev embedding $S^{1,2}\subset L^p$, \\
  \hspace{3,2cm} $\frac{1}{q}=\frac{1}{2}-\frac{1}{m}$ & \hspace{3,2cm}
  $\frac{1}{p}=\frac{1}{2}-\frac{1}{2n+2}$\\
   Conformal change $\wt{g}=u^{q-2}g$ & Change of contact form
  $\wt{\theta}=u^{p-2}\theta$\\
 \end{tabular}

\begin{tabular}{ll}
  \vspace{0,3cm}\underline{Conformal geometry} & \underline{CR geometry} \\
  Conformal invariant $\mu(M)$\hspace{0.8cm} &  CR invariant $\lambda(M)$ \\
  Yamabe equation:  & CR Yamabe equation:  \\
  \hspace{1cm} $a_m\Delta u+Ku=\mu u^{q-1}$ & \hspace{1cm} $b_n\Delta_b u+Ru=\lambda u^{p-1}$
\end{tabular}
\vspace{0,3cm}

Many authors have already expanded the previous list with important
contributions, as Gover and Graham did in \cite{GovGra}, where they
derived the CR analogues on CR manifolds of the GJMS operators
defined on Riemannian manifolds. For the original result on
Riemannian manifolds, see the paper by Graham, Jenne, Mason and
Sparling \cite{GJMS}.

The point of this paper is to provide a new characterization of
fully nonlinear CR invariant equations of the second order on the
Heisenberg group, thus adding another interesting similarity between
CR geometry and conformal geometry. We then also provide comparison
principles for solutions of families of fully nonlinear second order
operators on $\h^n$, which have suitable invariances.

The original result on the Euclidean space $\R^n$ was proved by A.
Li and the first author in \cite{LiLi}. There, among many other
results, they showed that any fully nonlinear conformally invariant
equation on $\R^n$ takes the form $$F\big(x,u,\nabla
u,\nabla^2u\big)=F\Big(0,1,0,-\frac{n-2}{2}A^u\Big),$$ where
\begin{equation}\label{55}
\quad
A^u:=-\frac{2}{n-2}u^{-\frac{n+2}{n-2}}\nabla^2u+\frac{2n}{(n-2)^2}
u^{-\frac{2n}{n-2}}\nabla u\otimes\nabla
u-\frac{2}{(n-2)^2}u^{-\frac{2n}{n-2}}|\nabla u|^2I_n,
\end{equation}
and $F(0,1,0,\cdot)$ is invariant under orthogonal conjugation, i.e.
$$F\Big(0,1,0,-\frac{n-2}{2}O^{-1}AO\Big)=F\Big(0,1,0,-\frac{n-2}{2}A\Big)$$
 for every real symmetric $n\times n$ matrix $A$ and real every orthogonal $n\times n$ matrix $O$.

The tensor $A^u$ is very closely related to the Schouten tensor of a
Riemannian manifold $(M,g)$, which is defined by setting
\begin{equation}\label{110}
A_g=\frac{1}{n-2}\bigg(\text{Ric}_g-\frac{R_g}{2(n-1)}g\bigg),
\end{equation}
where $\text{Ric}_g$ and $R_g$ denote the Ricci tensor and the
scalar curvature associated with $g$, respectively. Indeed, let $g_1
=u^\frac{4}{n-2}g$ be a conformal change of metrics on $M$; then, as
one can see in \cite{Via}, $$
A_{g_1}=-\frac{2}{n-2}u^{-1}\nabla^2_{g}u +
\frac{2n}{(n-2)^2}u^{-2}\nabla_gu\otimes\nabla_gu-\frac{2}{(n-2)^2}u^{-2}|\nabla_gu|^2_gg+A_g.$$
If one lets $g=u^\frac{4}{n-2}g_{\text{flat}}$, where
$g_{\text{flat}}$ denotes the Euclidean metric on $\R^n$, then by
the above transformation formula
$$A_g=u^\frac{4}{n-2}A^u_{ij}dx_idx_j,$$ where $A^u$ is given by \eqref{55}.

Lee derived in \cite{Lee} the analogue transformation law for the CR
Schouten tensor under a CR conformal change of the contact form
$\theta$ on a CR manifold.

\begin{rem}
Let $\N$ denote the set of positive integers. For any $N\in\N$ we
will denote by $I_N$ and $0_N$ the identity $N\times N$ matrix and
the zero $N\times N$ matrix respectively. We will denote by
$\text{Mat}(N,\R)$ the set of $N\times N$ real matrices and by
$\mathcal{S}^{N\times N}$ the set of real symmetric $N\times N$
matrices.

If $v,w\in\R^N$ for some $N\in\N$, we will denote by $v\otimes w$
the $N\times N$ real matrix
$$v\otimes w:=[v_iw_j]_{i,j=1,\ldots,N}.$$ Notice that for all
$v,w\in\R^N$ and all $A,B\in\text{Mat}(N,\R)$ one has
$$A\big(v\otimes w\big)B=\big(\big(Av\big)\otimes w\big)B=\big(Av\big)\otimes\big(B^Tw\big).$$
With some abuse of notation then we will simply write it as
$Av\otimes wB$.
\end{rem}

\subsection{\textbf{Fully nonlinear CR invariant equations of the second order on
$\h^n$}}

The Heisenberg group $\h^n$ is the set $\R^n\times\R^n\times\R$
endowed with the group action $\circ$ defined by
\begin{equation}\label{operation}
\xi\circ\hat{\xi}:=\Big(x+\hat{x},\,y+\hat{y},\,t+\hat{t}+2\sum_{i=1}^nx_i\hat{y}_i-y_i\hat{x}_i\Big)
\end{equation}
for any $\xi=(x,y,t)$, $\hat{\xi}=(\hat{x},\hat{y},\hat{t})$ in
$\h^n$, with $x=(x_1,\dots,x_n)$,
$\hat{x}=(\hat{x}_1,\dots,\hat{x}_n)$, $y=(y_1,\dots,y_n)$ and
$\hat{y}=(\hat{y}_1,\dots,\hat{y}_n)$ denoting elements of $\R^n$.
We will also use the notation $\xi=(z,t)$ with $z=x+iy$,
$z\in\C^n\simeq\R^n\times\R^n$. Let $Q:=2n+2$ denote the homogenous
dimension of $\h^n$, see also \cite{FolSte}. We consider the norm on
$\h^n$ defined by
\begin{equation}\label{norm}
\Hn{\xi}:=\Big[\Big(\sum_{i=1}^nx^2_i+y^2_i\Big)^2+t^2\Big]^\frac{1}{4}=\Big(|z|^4+t^2\Big)^\frac{1}{4}.
\end{equation}
The corresponding distance on $\h^n$ is defined accordingly by
setting
\begin{equation*}
d_{\!_H}(\xi,\hat{\xi}):=|\hat{\xi}^{-1}\circ\xi|_{\!_H}
\end{equation*}
where $\hat{\xi}^{-1}$ is the inverse of $\hat{\xi}$ with respect to
$\circ$, i.e. $\hat{\xi}^{-1}=-\hat{\xi}$. For every $\xi\in\h^n$
and $R>0$ we will use the notation
\begin{eqnarray}\label{ball}
D_R(\xi)&:=&\{\eta\in\h^n\,|\,d_{\!_H}(\xi,\eta)<R\}.
\end{eqnarray}

For any fixed $\hat{\xi}\in\h^n$ we will denote by
$\tau_{\hat{\xi}}:\h^n\rightarrow\h^n$ the \textbf{left translation}
on $\h^n$ by $\hat{\xi}$, defined by
\begin{equation}\label{trasl}
\tau_{\hat{\xi}}(\xi)=\hat{\xi}\circ\xi,
\end{equation}
where $\circ$ denotes the group action defined in \eqref{operation},
while for any $\lambda>0$ we will denote by
$\delta_\lambda:\h^n\rightarrow\h^n$ the \textbf{dilation} defined
by
\begin{equation}\label{dil}
\delta_\lambda(\xi):=(\lambda x,\, \lambda y,\, \lambda^2t),
\end{equation}
which satisfies
$$\delta_\lambda(\hat{\xi}\circ\xi)=\delta_\lambda(\hat{\xi})\circ\delta_\lambda(\xi)$$
for every $\xi$, $\hat{\xi}\in\h^n$ and every $\lambda>0$.

Notice that the norm on $\h^n$ defined by \eqref{norm} is
homogeneous of degree $1$ with respect to the dilations
$\delta_\lambda$, i.e.
$$|\delta_\lambda(\xi)|_{\!_H}=\lambda\Hn{\xi}\qquad\forall\,\xi\in\h^n,\,\lambda>0.$$

Another group of automorphisms of $\h^n$ is given by the action of
the n--dimensional unitary group $\mathcal{U}(n)$. Using the complex
numbers notation, its action is given by
\begin{equation}\label{rot}
\rho_{\!_M}(\xi)=\rho_{\!_M}(z,t):=(Mz,t)
\end{equation}
for any $M\in \mathcal{U}(n)$ and every $\xi=(z,t)\in\h^n$.

It is a known fact that a complex matrix $M\in\text{Mat}(n,\C)$
belongs to $\mathcal{U}(n)$, i.e. it satisfies
$M~\!\cdot\!~(\overline{M})^T=I_n$, if and only if the block matrix
$\wt{M}\in\text{Mat}(2n,\R)$ defined as in Theorem \ref{1} by
$$\wt{M}:=\left(
               \begin{array}{cc}
                 B & -C \\
                 C & B \\
               \end{array}
             \right),\qquad\text{with
             }\,B:=\text{Re}M,\,\,C:=\text{Im}M\,\in\,\text{Mat}(n,\R),
$$ belongs to $\mathcal{O}(2n)$,
i.e. one has $\wt{M}~\!\cdot\!~(\wt{M})^T=I_{2n}$. Using the real
numbers notation, one has
\begin{equation*}
\rho_{\!_M}(\xi)=\rho_{\!_M}(x,y,t)=(Bx-Cy,\,By+Cx,\,t)
\end{equation*}
for any $M\in \mathcal{U}(n)$ and every $\xi=(x,y,t)\in\h^n$. In the
case of the Euclidean space $\R^n$, the analogues of these maps are
the usual \textbf{rotations}, given by the action of the group
$\mathcal{O}(n)$ on $\R^n$.

We finally introduce the \textbf{inversion} map
$\iota:\h^n\rightarrow\h^n$ defined by
\begin{equation}\label{invers}
\iota(\xi)=\iota(x,y,t):=(x,-y,-t)
\end{equation}
for every $\xi=(x,y,t)\in\h^n$, and the map
$\varphi:\h^n\rightarrow\h^n$ defined by Jerison and Lee in
\cite{JerLee1} which we shall refer to as the \textbf{CR inversion}
and which is defined by the following relations:
\begin{equation}\label{24}
\varphi(\xi):=\wt\xi
\end{equation}
where $\wt\xi=(\wt{x},\wt{y},\wt{t})$ and
\begin{equation}\label{CR}
\wt{x}:=\frac{xt+y|z|^2}{\Hn{\xi}^4},\qquad\wt{y}:=\frac{yt-x|z|^2}{\Hn{\xi}^4},\qquad\wt{t}:=\frac{-t}{\Hn{\xi}^4}.
\end{equation}
We explicitly remark that
$\big|\varphi(\xi)\big|_{\!_H}=\frac{1}{\Hn\xi}$. The CR inversion
of $\h^n$ plays the role of the usual inversion with respect to the
unitary sphere in $\R^n$.

The elements of the group of automorphisms of $\h^n$ generated by
the left translations \eqref{trasl}, by the dilations \eqref{dil},
by the rotations \eqref{rot}, by the inversion map \eqref{invers}
and by the CR inversion \eqref{CR} are called \textbf{CR maps} on
$\h^n$.

For further references on these maps, on their definitions and their
properties we also refer to the works of Koranyi \cite{Kor}, Jerison
and Lee \cite{JerLee1} and Birindelli and Prajapat \cite{BirPra1},
\cite{BirPra2}.

The vector fields
\begin{equation}\label{vectorfields}
\begin{array}{rcl}
X_j&:=&\frac{\partial}{\partial x_j}+2y_j\frac{\partial}{\partial
t},\qquad j=1,\dots,n\\
Y_j&:=&\frac{\partial}{\partial y_j}-2x_j\frac{\partial}{\partial
t},\qquad j=1,\dots,n\\
T&:=&\frac{\partial}{\partial t}
\end{array}
\end{equation}
form a base of the Lie algebra of vector fields on the Heisenberg
group which are left invariant with respect to the group action
$\circ$. The Heisenberg gradient, or horizontal gradient, of a
regular function $u:\h^n\rightarrow\h^n$ is then defined by
$$\Hgrad{u}:=(X_1u,\dots,X_nu,Y_1u,\dots,Y_nu),$$ while its
Heisenberg hessian matrix is
\begin{equation}\label{Heishessian}
\begin{array}{rcl}
  \Hgrad^2{u}&:=&\left(
                 \begin{array}{ccc|ccc}
                   X_1X_1u & \cdots & X_nX_1u & Y_1X_1u & \cdots & Y_nX_1u \\
                   \vdots & \ddots & \vdots & \vdots & \ddots & \vdots \\
                   X_1X_nu & \cdots & X_nX_nu & Y_1X_nu & \cdots & Y_nX_nu \\
                   \hline X_1Y_1u & \cdots & X_nY_1u & Y_1Y_1u & \cdots & Y_nY_1u \\
                   \vdots & \ddots & \vdots & \vdots & \ddots & \vdots \\
                   X_1Y_nu & \cdots & X_nY_nu & Y_1Y_nu & \cdots & Y_nY_nu \\
                 \end{array}
               \right)\\
  &=&\left(\begin{array}{c|c}
                   X_jX_iu & Y_jX_iu \\
                   \hline X_jY_iu & Y_jY_iu
               \end{array}\right)_{i,j=1,\dots,n}
\end{array}
\end{equation}
We also define
\begin{equation}\label{Hess_symm}
\Hhesss{u(\xi)}:=\frac{1}{2}\left[\Hgrad^2u(\xi)+\left(\Hgrad^2u(\xi)\right)^T\right]
\end{equation}
which is the symmetric part of the matrix $\Hgrad^2u(\xi)$.

Now let $n\in\N$ and define $G,J\in\text{Mat}(2n,\R)$ by setting
\begin{equation}\label{45}
G:=\left(\begin{array}{cc}
               I_n & 0_n \\
               0_n & -I_n \\
         \end{array}\right),\qquad\qquad
J:=\left(\begin{array}{cc}
               0_n & I_n \\
               -I_n & 0_n \\
         \end{array}\right).
\end{equation}
We notice that $\Hgrad^2{u}\in\mathcal{S}^{2n\times2n}\oplus J\R$, see also section \ref{sec2_1}.

We now introduce the definition of CR invariance for an operator $F$
depending on $\xi,\,u,\,\Hgrad u,\,\Hgrad^2u$. We refer to section
\ref{notation} for some further basic facts concerning the
Heisenberg group.

\begin{defn}\label{def_conf_inv} Let $F\in\mathcal{C}^0\big(\h^n\times\R^+\times\R^{2n}\times
(\mathcal{S}^{2n\times2n}\oplus J\R)\big)$. Then $F(\cdot,u,\Hgrad
u,\Hgrad^2u)$ is CR invariant on the Heisenberg group $\h^n$ if for
any positive function $u\in\mathcal{C}^2(\h^n)$ and any CR map
$\psi:\h^n\rightarrow\h^n$ one has that
\begin{equation}\label{conf_inv}
F\big(\xi,u_\psi(\xi),\Hgrad u_\psi(\xi),\Hgrad^2u_\psi(\xi)\big)=
F\big(\psi(\xi),u(\psi(\xi)),\Hgrad
u(\psi(\xi)),\Hgrad^2u(\psi(\xi))\big),
\end{equation}
for every $\xi$, where the function $u_\psi$ is the transformed
function of $u$ through the CR map $\psi$, which is defined by
\begin{equation}\label{49}
u_\psi(\xi):=|J_\psi(\xi)|^\frac{Q-2}{2Q}u\big(\psi(\xi)\big),\qquad\xi\in\h^n.
\end{equation}
\end{defn}

The main results of the present paper are now contained in the
following theorems.

\begin{thm}\label{1}  Let $F\in\mathcal{C}^0\big(\h^n\times\R^+\times\R^{2n}\times
(\mathcal{S}^{2n\times2n}\oplus J\R)\big)$. Then $F(\cdot,u,\Hgrad
u,\Hgrad^2u)$ is CR invariant on the Heisenberg group $\h^n$ if and
only if
\begin{equation}\label{68}
F\big(\cdot,u,\Hgrad u,\Hgrad^2u\big)\equiv
F\Big(0,1,0,-\frac{Q-2}{2}A^u\Big),
\end{equation}
where
\begin{eqnarray}
\label{A^u}\hspace{1cm}A^u\!&\!:=\!&\!-\frac{2}{Q-2}u^{-\frac{Q+2}{Q-2}}\Hhesss{u}+\frac{2Q}{(Q-2)^2}u^{-\frac{2Q}{Q-2}}\Hgrad
    u\otimes\Hgrad u\\
\nonumber\!&\!\!&\!\hspace{1cm}-\frac{4}{(Q-2)^2}u^{-\frac{2Q}{Q-2}}J\Hgrad
    u\otimes J\Hgrad u-\frac{2}{(Q-2)^2}u^{-\frac{2Q}{Q-2}}|\Hgrad
    u|^2I_{2n},
\end{eqnarray}
and moreover for every $A\in\mathcal{S}^{2n\times2n}\oplus J\R$ one
has
\begin{itemize}
\item[i)] $F(0,1,0,A)=F(0,1,0,\wt{M}^{T}A\wt{M})$ for every
unitary matrix $M=B+iC\in\mathcal{U}(n)$, where we have set
$\wt{M}:=\left(\begin{array}{cc}
                   B & -C \\
                   C & B \\
                \end{array}\right),$
\item[ii)] $F(0,1,0,A)=F(0,1,0,GAG)$,
\item[iii)] $F(0,1,0,A)=F(0,1,0,A+\alpha J)$ for every
$\alpha\in\R$,
\end{itemize}
with $J$, $G$ being defined as in \eqref{45}.
\end{thm}

We want to stress here that $A^u\in\mathcal{S}^{2n\times2n}$, thus
it always has real eigenvalues, even if the Heisenberg hessian
matrix $\Hgrad^2 u$ in general is not symmetric. Let
$$\lambda(A^u)=(\lambda_1(A^u),\ldots,\lambda_{2n}(A^u))$$ denote the
eigenvalues of $A^u$. Using Theorem \ref{1} we can then provide some
examples of fully nonlinear CR invariant operators of the second
order on the Heisenberg group. Indeed, for $k=1,\ldots,2n$, let
$$\omega_k(\lambda)=\sum_{1\leq i_1<\ldots<i_k\leq 2n}
\lambda_{i_1}\cdots\lambda_{i_k},\qquad \lambda=
(\lambda_1,\ldots,\lambda_{2n})\in\R^{2n}$$ denote the
$k^{\text{th}}$ symmetric function on $\R^{2n}$. Then
$\omega_k\big(\lambda(A^u)\big)$ is a fully nonlinear CR invariant
differential operator of the second order on $\h^n$.

Similar operators, involving the tensors in equations \eqref{55} and
\eqref{110}, were studied in the context of conformal geometry on
$\R^n$ and on more general Riemannian manifolds by many authors, see
the works by Viaclovsky \cite{Via} and \cite{Via2}, the papers by
Chang, Gursky, and Yang \cite{ChaGurYan1} and \cite{ChaGurYan2} and
the works by A. Li and the first author \cite{LiLi} and
\cite{LiLi2}, and the references therein.

\vspace{0,2cm}
\subsection{\textbf{Comparison principle on a domain $\Omega\subset\h^n$ for fully nonlinear CR invariant equations}}

\begin{rem}
If $A,B\in\text{Mat}(N,\R)$ for some $N\in\N$, we will write $A\geq
B$ if
$$\crochet{\xi}{A\xi}_{\R^N}\geq\crochet{\xi}{B\xi}_{\R^N}\quad\text{for
every }\xi\in\R^N,$$ where we denoted with
$\crochet{\,\cdot}{\cdot}_{\R^N}$ the usual scalar product in
$\R^N$. If $A-B$ is a diagonalizable matrix, this is equivalent to
$A-B$ having nonnegative eigenvalues, i.e. to $A-B$ being
nonnegative definite.
\end{rem}

Let $\Sigma\subset\mathcal{S}^{2n\times2n}$ be an open set of
matrices such that
\begin{equation}\label{14}
\begin{array}{rl}
\text{i)}&A\in\overline{\Sigma},\,c\in\R^+\Longrightarrow
   cA\in\overline{\Sigma},\\
\text{ii)}&A\in\overline{\Sigma},\,B\in\mathcal{S}^{2n\times2n}\text{
and }B>0
   \Longrightarrow A+B\in\Sigma.
\end{array}
\end{equation}
Notice that condition ii) in particular implies
\begin{equation*}
\text{iii)}\,\,\,\,A\in\overline{\Sigma},\,B\in\mathcal{S}^{2n\times2n}\text{
and }B\geq0\Longrightarrow A+B\in\overline{\Sigma}.
\end{equation*}

\begin{thm}\label{13}
Let $\Omega\subset\h^n$ be a domain, let
$\Sigma\subset\mathcal{S}^{2n\times2n}$ be an open set of matrices
satisfying condition \eqref{14} and let
$u,w\in\mathcal{C}^2(\Omega)\cap\mathcal{C}^0(\overline{\Omega})$.
Assume that $u,w>0$ on $\overline{\Omega}$,
$A^u\in\overline{\Sigma}$ and $A^w\in\Sigma^c$ for every
$\xi\in\Omega$. Then
\begin{itemize}
\item[i)] if $u\geq w$ on
$\partial\Omega$, $u\geq w$ in $\overline{\Omega}$,
\item[ii)] if $u>w$ on
$\partial\Omega$, $u> w$ in $\overline{\Omega}$.
\end{itemize}
\end{thm}

\vspace{0,2cm}
\subsection{\textbf{Comparison principle for a family of fully nonlinear equations on a punctured ball
$\mathbf{D_R(\xi_0)\subset\h^n}$}}

\begin{rem} In this section we will consider an operator
$T\!\in\!\mathcal{C}^1\!\big(\R^+\!\times\!\R^{2n}\!\times\!(\mathcal{S}^{2n\times2n}\oplus
J\R)\!\big)$ satisfying the following assumptions
\begin{itemize}
\item[i)] $T=T(s,v,U)$ is elliptic with respect to the family of vector
fields $X_1,\ldots,X_n,\,$ $Y_1,\ldots,Y_n$, i.e.
\begin{equation}\label{106}
\left[-\frac{\partial T}{\partial
U_{ij}}(s,v,U)\right]_{i,j=1,\ldots,2n}>0\qquad\text{ on
}\R^+\!\times\R^{2n}\!\times (\mathcal{S}^{2n\times2n}\oplus J\R),
\end{equation}
that is for every $[a,b]\subset\R^+$, $K_1\subset\R^{2n}$ compact,
$K_2\subset\mathcal{S}^{2n\times2n}\oplus J\R$ compact there exists
$\beta>0$ such that $$\crochet{W}{-\frac{\partial T}{\partial
U_{ij}}(s,v,U)\cdot W}_{\R^{2n}}\geq\beta|W|^2$$ for every
$W\in\R^{2n}$, $s\in[a,b]$, $v\in K_1$, $U\in K_2$;
\item[ii)] $T$ is invariant with respect to dilations in $\h^n$, i.e. for every
$\lambda>0$ and every positive function $u\in\mathcal{C}^2(\h^n)$
one has
$$T\big(u_{\delta_\lambda}(\xi),\Hgrad u_{\delta_\lambda}(\xi),\Hgrad^2u_{\delta_\lambda}(\xi)\big)=
T\big(u\big(\delta_\lambda(\xi)\big),\Hgrad
u\big(\delta_\lambda(\xi)\big),\Hgrad^2u\big(\delta_\lambda(\xi)\big)\big)$$
for every $\xi\in\h^n$.
\end{itemize}
Notice that since the operator $T$ does not explicitly depend on
$\xi\in\h^n$, it is automatically invariant with respect to left
translations in $\h^n$, i.e. for every $\hat\xi\in\h^n$ and every
positive function $u\in\mathcal{C}^2(\h^n)$ one has
$$T\big(u_{\tau_{\hat\xi}}(\xi),\Hgrad u_{\tau_{\hat\xi}}(\xi),\Hgrad^2u_{\tau_{\hat\xi}}(\xi)\big)=
T\big(u\big(\tau_{\hat\xi}(\xi)\big),\Hgrad
u\big(\tau_{\hat\xi}(\xi)\big),\Hgrad^2u\big(\tau_{\hat\xi}(\xi)\big)\big)$$
for every $\xi\in\h^n$.
\end{rem}

\begin{thm}\label{94} Assume that the operator $T$ satisfies the hypotheses
above. Consider $D_2(\xi_0)\subset\h^n$ and let
$u\in\mathcal{C}^2(D_2(\xi_0)\setminus\{\xi_0\})$,
$w\in\mathcal{C}^2(D_2(\xi_0))$ be such that
\begin{itemize}
\item[i)] $u>w$ in $D_2(\xi_0)\setminus\{\xi_0\}$ and $w>0$ in
$D_2(\xi_0)$,
\item[ii)] $\Delta_{\!_H}u\leq0$ in $D_2(\xi_0)\setminus\{\xi_0\}$,
\item[iii)] $T(u,\Hgrad u,\Hgrad^2u)\geq0\geq T(w,\Hgrad
w,\Hgrad^2w)$ in $D_2(\xi_0)\setminus\{\xi_0\}$.
\end{itemize}
Then
$\displaystyle\liminf_{\xi\rightarrow\xi_0}\big(u(\xi)-w(\xi)\big)>0$.
\end{thm}
This result is an analogue in the Heisenberg group setting of the
original result proved in the Euclidean setting by the first author
(see Theorem 1.7 in \cite{YLi2}). From Theorem \ref{1} and Theorem
\ref{94}, the following Corollary immediately follows.

\begin{cor}\label{22}
Let $F(A^u)$ be a CR invariant operator on the Heisenberg group,
with $A^u$ being defined as in equation \eqref{A^u} and
$F\in\mathcal{C}^1(\mathcal{S}^{2n\times2n})$. Assume that for every
$A\in\mathcal{S}^{2n\times2n}$ one has
$$\left[\frac{\partial F}{\partial
A_{ij}}(A)\right]_{i,j=1,\ldots,2n}>0.$$ Consider
$D_2(\xi_0)\subset\h^n$ and let
$u\in\mathcal{C}^2(D_2(\xi_0)\setminus\{\xi_0\})$,
$v\in\mathcal{C}^2(D_2(\xi_0))$ be such that
\begin{itemize}
\item[i)] $u>w$ in $D_2(\xi_0)\setminus\{\xi_0\}$ and $w>0$ in
$D_2(\xi_0)$,
\item[ii)] $\Delta_{\!_H}u\leq0$ in $D_2(\xi_0)\setminus\{\xi_0\}$,
\item[iii)] $F(A^u)\geq0\geq F(A^w)$ in $D_2(\xi_0)\setminus\{\xi_0\}$.
\end{itemize}
Then
$\displaystyle\liminf_{\xi\rightarrow\xi_0}\big(u(\xi)-w(\xi)\big)>0$.
\end{cor}

\vspace{0,3cm}

The paper is organized as follows: in section \ref{notation} we
introduce notations, definitions and some known facts about the
Heisenberg group and the sublaplacian operator defined on it, which
we are going to use throughout the paper. Section \ref{notation}
also contains the formulae for $\Hgrad u_\psi$ and $\Hgrad^2
u_\psi$, when $u\in\mathcal{C}^2(\h^n)$ and $\psi$ is any of the
generators of the group of CR maps on $\h^n$.

In section \ref{47} we give the proof of Theorem \ref{1}, in section
\ref{56} we prove Theorem \ref{13} while in section \ref{95} we give
the proof of Theorem \ref{94}. Section \ref{56} also contains
another result of interest, where we study ``the first variation''
of the operator $A^u$ defined in \eqref{A^u}, and which is the
equivalent on the Heisenberg group of an analogous lemma proved in
the Euclidean setting by the first author (see lemma 3.7 in
\cite{YLi}).

Finally in section \ref{57} we collect some technical results, which
are used in the previous sections.

\section{\textbf{Notation and preliminary facts}}\label{notation}

For future use we notice here that if we denote by $J_\psi(\xi)$ the
Jacobian matrix of a CR map $\psi:\h^n\rightarrow\h^n$ evaluated at
$\xi\in\h^n$ and by $|J_\psi(\xi)|$ its determinant, then we have
\begin{eqnarray}
\nonumber&|J_{\tau_{\hat{\xi}}}(\xi)|=1,\qquad|J_{\rho_{\!_M}}(\xi)|=1,\qquad|J_\iota(\xi)|=1&\\
\nonumber&|J_{\delta_\lambda}(\xi)|=\lambda^Q,\qquad|J_\varphi(\xi)|=\frac{1}{\Hn{\xi}^{2Q}},&
\end{eqnarray}
where we recall that $Q=2n+2$ denote the homogenous dimension of
$\h^n$.

Next we recall that, for any CR map $\psi$ on $\h^n$ and any function
$u:\h^n\rightarrow\h^n$, the transformed function $u_\psi$ is defined as in \eqref{49} by
\begin{equation*}
u_\psi(\xi):=|J_\psi(\xi)|^\frac{Q-2}{2Q}u\big(\psi(\xi)\big),\qquad\xi\in\h^n.
\end{equation*}

Then we have
\begin{eqnarray}
\label{4}\hspace{1cm}
&u_{\tau_{\hat{\xi}}}(\xi)=u(\hat{\xi}\circ\xi),\quad
u_{\rho_{\!_M}}(\xi)=u(Bx-Cy,Cx+By,t),\quad
u_{\iota}(\xi)=u(x,-y,-t)&\\
\nonumber
\hspace{1cm}&u_{\delta_\lambda}(\xi)=\lambda^\frac{Q-2}{2}u(\lambda
x,\lambda y,\lambda^2t),\qquad
u_{\varphi}(\xi)=\frac{1}{\Hn{\xi}^{Q-2}}u(\wt{x},\wt{y},\wt{t}).&
\end{eqnarray}

\subsection{\textbf{The sublaplacian on the Heisenberg
group}}\label{sec2_1}

Consider the vector fields $X_j,Y_j$ for $j=1,\ldots,n$ defined in
\eqref{vectorfields}. The sublaplacian on the Heisenberg group is
the linear differential operator of the second order defined by
\begin{eqnarray}
\nonumber\Delta_{\!_H}u&:=&\sum_{j=1}^n X_j^2u+Y_j^2u\\
\nonumber          &=&\sum_{j=1}^n\frac{\partial^2 u}{\partial
                          x_j^2}+\frac{\partial^2 u}{\partial
                          y_j^2}+4y_j\frac{\partial^2 u}{\partial
                          x_j\partial t}-4x_j\frac{\partial^2 u}{\partial
                          y_j\partial t}+4(x_j^2+y^2_j)\frac{\partial^2 u}{\partial
                          t^2}.
\end{eqnarray}

The sublaplacian is the trace of the Heisenberg hessian matrix
defined in \eqref{Heishessian} and it is degenerate elliptic.
Furthermore it has divergence form. Indeed one has
$$\Delta_{\!_H}u=\text{div}(A(z)\nabla u),$$
where here $\nabla u$ denotes the gradient of $u$ in $\R^{2n+1}$ and
$$A(z)=A(x,y):=\left(
    \begin{array}{ccc}
      I_n & 0_n & 2y \\
      0_n & I_n & -2x \\
      2y & -2x & 4|z|^2 \\
    \end{array}
  \right),$$
with $I_n$ and $0_n$ denoting respectively the identity matrix and
the zero matrix in $\R^n$.

Now notice that the smooth vector fields $X_j$, $Y_j$, $j=1,\dots,n$
satisfy the following commutation relations
\begin{equation}\label{commutation}
\left[X_i,Y_j\right]=-4T\delta_{ij},\quad\left[X_i,X_j\right]=\left[Y_i,Y_j\right]=0\qquad\,\forall\,i,j=1,\dots,n
\end{equation}
and hence they and their first order commutators span the whole Lie
Algebra. It is then a consequence of the classical theorem of
H\"{o}rmander \cite{Hor} that the sublaplacian is hypoelliptic, i.e.
if $\Delta_{\!_H}u\in\mathcal{C}^\infty$ then
$u\in\mathcal{C}^\infty$. Moreover $\Delta_{\!_H}$ satisfies the
Strong Maximum Principle, as one finds in the work of Bony
\cite{Bony}.

\begin{rem}
Notice that, by the commutation relations \eqref{commutation}, the
Heisenberg hessian matrix of a regular function $u$ is not
symmetric, in general. Indeed, it's easy to see that
\begin{eqnarray}
\label{78}\Hgrad^2u(\xi)&=&\Hhesss{u(\xi)}+2Tu(\xi)J,
\end{eqnarray}
with $J\in\text{Mat}(2n,\R)$ being defined as in \eqref{45} by
\begin{equation*}
J:=\left(
    \begin{array}{cc}
      0_n & I_n \\
      -I_n & 0_n \\
    \end{array}
\right).
\end{equation*}
\end{rem}

\begin{rem}\label{30}
Using the definition of the matrix $J$ in \eqref{45}, for every
regular function $u$ one has
\begin{eqnarray}
\nonumber\Hgrad{u}\!&\!=\!&\!\nabla_zu+2\frac{\partial u}{\partial t}Jz\\
\nonumber\Hgrad^2u\!&\!=\!&\!\nabla^2_zu+2\frac{\partial}{\partial
  t}\nabla_zu\otimes Jz+2Jz\otimes\frac{\partial}{\partial t}\nabla_zu+
  4\frac{\partial^2u}{\partial t^2}Jz\otimes Jz+2\frac{\partial u}{\partial t}J\\
\nonumber\Hhesss u\!&\!=\!&\!\nabla^2_zu+2\frac{\partial}{\partial
  t}\nabla_zu\otimes Jz+2Jz\otimes\frac{\partial}{\partial t}\nabla_zu+
  4\frac{\partial^2u}{\partial t^2}Jz\otimes Jz\\
\nonumber \Delta_{\!_H}u\!&\!=\!&\!\Delta_zu+4|z|^2
  \frac{\partial^2u}{\partial t^2}+ 4\frac{\partial}{\partial
  t}\crochet{Jz}{\nabla_zu}_{\R^{2n}},
\end{eqnarray}
where $\Delta_z$, $\nabla^2_z$ and $\nabla_z$ are respectively the
ordinary Laplace operator, the Hessian matrix and the gradient with
respect to the variables $z=(x,y)\in\R^{2n}$. This in particular
implies
$$\Hgrad{u}(0)=\nabla_zu(0),\qquad\Hgrad^2u(0)=\nabla^2_zu(0)+2\frac{\partial u}{\partial t}(0)J.$$
\end{rem}

\begin{rem}
It is not difficult to see that
\begin{eqnarray}
\nonumber\Hgrad u_{\tau_{\hat{\xi}}}(\xi)\!&\!=\!&\!\Hgrad u(\hat{\xi}\circ\xi)\\
\label{grad_tr}\Hgrad u_{\rho_{\!_M}}(\xi)\!&\!=\!&\!\wt{M}^T\cdot\Hgrad u(Bx-Cy,By+Cx,t)\\
\nonumber\Hgrad u_{\iota}(\xi)\!&\!=\!&\!G\cdot\Hgrad u(x,-y,-t)\\
\nonumber\Hgrad
u_{\delta_\lambda}(\xi)\!&\!=\!&\!\lambda^\frac{Q}{2}\Hgrad
u(\lambda x,\lambda y,\lambda^2t),
\end{eqnarray}
with $G\in\text{Mat}(2n,\R)$ being defined as in Theorem \ref{1} by
$G=\left(
      \begin{array}{cc}
        I_n & 0_n \\
        0_n & -I_n \\
      \end{array}
    \right),$
and that
\begin{eqnarray}
\nonumber\Hgrad^2 u_{\tau_{\hat{\xi}}}(\xi)\!&\!=\!&\!\Hgrad^2 u(\hat{\xi}\circ\xi)\\
\label{Hess_tr}\Hgrad^2 u_{\rho_{\!_M}}(\xi)\!&\!=\!&\!\wt{M}^T\cdot\Hgrad^2 u(Bx-Cy,By+Cx,t)\cdot\wt{M}\\
\nonumber\Hgrad^2 u_{\iota}(\xi)\!&\!=\!&\!G\cdot\Hgrad^2 u(x,-y,-t)\cdot G\\
\nonumber\Hgrad^2
u_{\delta_\lambda}(\xi)\!&\!=\!&\!\lambda^\frac{Q+2}{2}\Hgrad^2
u(\lambda x,\lambda y,\lambda^2t).
\end{eqnarray}

Then it follows that
\begin{eqnarray}
\nonumber\Delta_{\!_H} u_{\tau_{\hat{\xi}}}(\xi)\!&\!=\!&\!\Delta_{\!_H} u(\hat{\xi}\circ\xi)\\
\label{Lapl_tr}\Delta_{\!_H} u_{\rho_{\!_M}}(\xi)\!&\!=\!&\!\Delta_{\!_H} u(Bx-Cy,By+Cx,t)\\
\nonumber\Delta_{\!_H} u_{\iota}(\xi)\!&\!=\!&\!\Delta_{\!_H} u(x,-y,-t)\\
\nonumber\Delta_{\!_H}u_{\delta_\lambda}(\xi)\!&\!=\!&\!\lambda^\frac{Q+2}{2}\Delta_{\!_H}
u(\lambda x,\lambda y,\lambda^2t).
\end{eqnarray}
\end{rem}

\subsection{\textbf{A useful CR map}}\label{sec2_2}

Instead of using the CR inversion $\varphi$ defined in \eqref{24} as
one of the generators of the group of CR maps on $\h^n$, throughout
the rest of the paper we will use the map
$\check{\varphi}:=\varphi\circ\iota$, i.e.
$\check{\varphi}(\xi)=\big(\check{x},\,\check{y},\,\check{t}\big)$
for every $\xi\in\h^n$ with $(\check{x},\,\check{y},\,\check{t})$
being in turn defined by
\begin{equation}\label{48}
\check{x}=-\frac{xt+y|z|^2}{\Hn{\xi}^4},\quad\check{y}=\frac{yt-x|z|^2}{\Hn{\xi}^4},\quad\check{t}=\frac{t}{\Hn{\xi}^4}.
\end{equation}
We make this choice because
$\check\varphi\big(\check\varphi(\xi)\big)=\xi$, while
$\varphi\big(\varphi(\xi)\big)=(-x,-y,t)$ for every
$\xi=(x,y,t)\in\h^n\setminus\{0\}$.

We now derive the transformation formulae for $\Hgrad
u_{\check\varphi}$ and $\Hgrad^2u_{\check\varphi}$.

\begin{pro}\label{51}
Let $u\in\mathcal{C}^2(\h^n)$, then
$u_{\check{\varphi}}(\xi)=\frac{1}{\Hn{\xi}^{Q-2}}u(\check{x},\check{y},\check{t})$
for every $\xi\in\h^n\setminus\{0\}$. Moreover for every
$j=1,\ldots,n$ and every $\xi=(x,y,t)\in\h^n\setminus\{0\}$ one has
\begin{eqnarray}
\nonumber
X_ju_{\check{\varphi}}(\xi)\!\!&\!\!=\!\!&\!\!-(Q-2)\Hn{\xi}^{-(Q+2)}u\big(\check{\varphi}(\xi)\big)\big(|z|^2x_j+ty_j\big)\\
\nonumber
\!\!&\!\!\!\!&\!\!\hspace{0,3cm}+\Hn{\xi}^{-(Q+6)}\frac{\partial
u}{\partial
t}(\check{\varphi}\big(\xi)\big)\big(2\big(|z|^4-t^2\big)y_j-4t|z|^2x_j\big)\\
\label{60}\!\!&\!\!\!\!&\!\!\hspace{0,3cm}+\Hn{\xi}^{-(Q+6)}\Big[\sum_{h=1}^n\frac{\partial
u}{\partial
x_h}\big(\check{\varphi}(\xi)\big)\Big(-\delta_{\!jh}t\Hn{\xi}^4+2\big(|z|^4-t^2\big)(x_jy_h-y_jx_h)\\
\nonumber\!\!&\!\!\!\!&\!\!\hspace{7,9cm}+4t|z|^2(x_jx_h+y_jy_h)\Big)\Big]\\
\nonumber\!\!&\!\!\!\!&\!\!\hspace{0,3cm}+\Hn{\xi}^{-(Q+6)}\Big[\sum_{h=1}^n\frac{\partial
u}{\partial
y_h}\big(\check{\varphi}(\xi)\big)\Big(-\delta_{\!jh}|z|^2\Hn{\xi}^4+2\big(|z|^4-t^2\big)(x_jx_h+y_jy_h)\\
\nonumber\!\!&\!\!\!\!&\!\!\hspace{7,9cm}+4t|z|^2(y_jx_h-x_jy_h)\Big)\Big]\\
\nonumber
Y_ju_{\check{\varphi}}(\xi)\!\!&\!\!=\!\!&\!\!-(Q-2)\Hn{\xi}^{-(Q+2)}u\big(\check{\varphi}(\xi)\big)\big(|z|^2y_j-tx_j\big)\\
\nonumber
\!\!&\!\!\!\!&\!\!\hspace{0,3cm}+\Hn{\xi}^{-(Q+6)}\frac{\partial
u}{\partial
t}(\check{\varphi}\big(\xi)\big)\big(2\big(t^2-|z|^4\big)x_j-4t|z|^2y_j\big)\\
\label{61}\!\!&\!\!\!\!&\!\!\hspace{0,3cm}+\Hn{\xi}^{-(Q+6)}\Big[\sum_{h=1}^n\frac{\partial
u}{\partial
x_h}\big(\check{\varphi}(\xi)\big)\Big(-\delta_{\!jh}|z|^2\Hn{\xi}^4+2\big(|z|^4-t^2\big)(y_jy_h+x_jx_h)\\
\nonumber\!\!&\!\!\!\!&\!\!\hspace{7,9cm}+4t|z|^2(y_jx_h-x_jy_h)\Big)\Big]
\end{eqnarray}
\begin{eqnarray}
\nonumber\hspace{2cm}\!\!&\!\!\!\!&\!\!\hspace{0,3cm}+\Hn{\xi}^{-(Q+6)}\Big[\sum_{h=1}^n\frac{\partial
u}{\partial
y_h}\big(\check{\varphi}(\xi)\big)\Big(\delta_{\!jh}t\Hn{\xi}^4+2\big(|z|^4-t^2\big)(y_jx_h-x_jy_h)\\
\nonumber\!\!&\!\!\!\!&\!\!\hspace{7,9cm}-4t|z|^2(x_jx_h+y_jy_h)\Big)\Big].
\end{eqnarray}
\end{pro}
\textbf{Proof:} By definition \eqref{49} one has
$u_{\check{\varphi}}(\xi)=|J_{\check{\varphi}}(\xi)|^\frac{Q-2}{2Q}u(\check{\varphi}(\xi))$.
Since for every $\xi\in\h^n$ one has that
$$|J_{\check{\varphi}}(\xi)|=|J_{\varphi\circ\iota}(\xi)|=|J_\varphi(\xi)||J_\iota(\xi)|=|J_\varphi(\xi)|=\frac{1}{\Hn{\xi}^{2Q}},$$
we get
$u_{\check{\varphi}}(\xi)=\frac{1}{\Hn{\xi}^{Q-2}}u(\check{x},\check{y},\check{t})$
as claimed.

Then one also gets
\begin{eqnarray}
\nonumber
X_ju_{\check{\varphi}}(\xi)\!\!&\!\!=\!\!&\!\!u(\check{x},\check{y},\check{t})X_j\bigg(\frac{1}{\Hn{\xi}^{Q-2}}\bigg)+
\frac{1}{\Hn{\xi}^{Q-2}}X_j\Big(u(\check{x},\check{y},\check{t})\Big)\\
\nonumber\!\!&\!\!=\!\!&\!\!-(Q-2)\Hn{\xi}^{-(Q+2)}u(\check{x},\check{y},\check{t})\big(|z|^2x_j+ty_j\big)+\\
\nonumber\!\!&\!\!\!\!&\!\!\frac{1}{\Hn{\xi}^{Q-2}}\Big(\sum_{h=1}^n\frac{\partial
u}{\partial
x_h}(\check{x},\check{y},\check{t})X_j(\check{x}_h)+\sum_{h=1}^n\frac{\partial
u}{\partial
y_h}(\check{x},\check{y},\check{t})X_j(\check{y}_h)+\frac{\partial
u}{\partial t}(\check{x},\check{y},\check{t})X_j(\check{t})\Big),
\end{eqnarray}
and in a similar way
\begin{eqnarray}
\nonumber
Y_ju_{\check{\varphi}}(\xi)\!&\!=\!&\!u(\check{x},\check{y},\check{t})Y_j\bigg(\frac{1}{\Hn{\xi}^{Q-2}}\bigg)+
\frac{1}{\Hn{\xi}^{Q-2}}Y_j\Big(u(\check{x},\check{y},\check{t})\Big)\\
\nonumber\!&\!=\!&\!-(Q-2)\Hn{\xi}^{-(Q+2)}u(\check{x},\check{y},\check{t})\big(|z|^2y_j-tx_j\big)+\\
\nonumber\!&\!\!&\!\frac{1}{\Hn{\xi}^{Q-2}}\Big(\sum_{h=1}^n\frac{\partial
u}{\partial
x_h}(\check{x},\check{y},\check{t})Y_j(\check{x}_h)+\sum_{h=1}^n\frac{\partial
u}{\partial
y_h}(\check{x},\check{y},\check{t})Y_j(\check{y}_h)+\frac{\partial
u}{\partial t}(\check{x},\check{y},\check{t})Y_j(\check{t})\Big).
\end{eqnarray}
Now it suffices to use the formulae for
$X_j(\check{x}_h),\,X_j(\check{y}_h),\,X_j(\check{t}),\,Y_j(\check{x}_h),\,Y_j(\check{y}_h)$
and $Y_j(\check{t})$, which are provided in lemma \ref{50} in the
appendix, in order to conclude the proof.$\qquad\Box$

\begin{rem}\label{59} We can rewrite formulae \eqref{60} and
\eqref{61} in the following form
\begin{equation}\label{62}
\Hgrad{u_{\check{\varphi}}(\xi)}=-\frac{(Q-2)}{\Hn{\xi}^{Q+2}}u\big(\check{\varphi}(\xi)\big)\big(|z|^2z+tJz\big)+
\frac{1}{\Hn{\xi}^Q}\,E\cdot\Hgrad{u}\big(\check{\varphi}(\xi)\big),
\end{equation}
where $J\in\text{Mat}(2n,\R)$ is the matrix defined in \eqref{45}
and where $E\in\text{Mat}(2n,\R)$ is the block matrix defined by
\begin{equation}\label{63}
E:=\Hn{\xi}^2\left(\begin{array}{c|c}
                   X_j(\check{x}_h) & X_j(\check{y}_h) \\
                   \hline Y_j(\check{x}_h) & Y_j(\check{y}_h)
          \end{array}\right)_{j,h=1,...,n}.
\end{equation}
Using lemma \ref{50} we see that the matrix $E$ can be written in
the form
\begin{eqnarray}
\label{64}
E&=&\bigg(-\frac{t}{\Hn{\xi}^2}I_{2n}+\frac{|z|^2}{\Hn{\xi}^2}J+\frac{1}{\Hn{\xi}^6}\Big(2(|z|^4-t^2)(z\otimes
Jz-Jz\otimes z)+\\
\nonumber&&\hspace{6cm}4t|z|^2(z\otimes z+Jz\otimes Jz)\Big)\bigg)\cdot G\\
\nonumber&=& \left(\begin{array}{cc}
                   R & S \\
                   S & -R \\
                \end{array}\right)\,\,\,=\,\,\,\left(\begin{array}{cc}
                   R & -S \\
                   S & R \\
                \end{array}\right)\cdot G\,\,\,=\,\,\,G\cdot\left(\begin{array}{cc}
                   R & S \\
                   -S & R \\
                \end{array}\right),
\end{eqnarray}
where $G\in\text{Mat}(2n,\R)$ is the matrix defined in \eqref{45}
and where $R,S\in\text{Mat}(n,\R)$ are given by
\begin{eqnarray}
\nonumber
R&:=&\Hn{\xi}^2\Big(X_j(\check{x}_h)\Big)_{j,h=1,...,n}\,\,=\,\,-\Hn{\xi}^2\Big(Y_j(\check{y}_h)\Big)_{j,h=1,...,n}\\
\nonumber
&=&-\frac{t}{\Hn{\xi}^2}I_n+\frac{1}{\Hn{\xi}^6}\Big(2(|z|^4-t^2)(x\otimes
y-y\otimes x)+4t|z|^2(x\otimes x+y\otimes y)\Big)\\
\nonumber
S&:=&\Hn{\xi}^2\Big(X_j(\check{y}_h)\Big)_{j,h=1,...,n}\,\,=\,\,\Hn{\xi}^2\Big(Y_j(\check{x}_h)\Big)_{j,h=1,...,n}\\
\nonumber
&=&-\frac{|z|^2}{\Hn{\xi}^2}I_n+\frac{1}{\Hn{\xi}^6}\Big(2(|z|^4-t^2)(x\otimes
x+y\otimes y)+4t|z|^2(y\otimes x-x\otimes y)\Big).
\end{eqnarray}
With these definitions it's easy to see that $R+iS\in\mathcal{U}(n)$
and hence that $E\in\mathcal{O}(2n)$.
\end{rem}

\begin{cor}\label{66}
Let $u\in\mathcal{C}^2(\h^n)$, then
\begin{eqnarray}
\nonumber\Hgrad^2u_{\check\varphi}(\xi)&=&\frac{Q^2-4}{\Hn{\xi}^{Q+6}}u\big(\check\varphi(\xi)\big)\big(|z|^2z+tJz\big)\otimes\big(|z|^2z+tJz\big)\\
\nonumber&&-\frac{(Q-2)}{\Hn{\xi}^{Q+2}}u\big(\check{\varphi}(\xi)\big)\big(|z|^2I_{2n}+tJ+2z\otimes
z+2Jz\otimes
Jz\big)\\
\label{67}&&-\frac{(Q-2)}{\Hn{\xi}^{Q+4}}\big(|z|^2z+tJz\big)\otimes
\big(E\cdot\Hgrad u\big(\check{\varphi}(\xi)\big)\big)\\
\nonumber&&-\frac{(Q-2)}{\Hn{\xi}^{Q+4}}\big(E\cdot\Hgrad
u\big(\check{\varphi}(\xi)\big)\big)\otimes\big(|z|^2z+tJz\big)\\
\nonumber&&+\frac{1}{\Hn{\xi}^{Q+2}}E\cdot\Hgrad^2u\big(\check\varphi(\xi)\big)\cdot
E^T\\
\nonumber&&+\frac{1}{\Hn{\xi}^{Q-2}}\sum_{h=1}^n\Hgrad^2(\check{x}_h)X_hu\big(\check\varphi(\xi)\big)+\Hgrad^2(\check{y}_h)Y_hu\big(\check\varphi(\xi)\big)
\end{eqnarray}
for every $\xi\in\h^n\setminus\{0\}$, where $E$ is the matrix
defined in \eqref{63}.
\end{cor}

\textbf{Proof:} In order to obtain formula \eqref{67} it's
sufficient to recall the definition of the Heisenberg hessian matrix
$\Hgrad^2$ given in formula \eqref{Heishessian} and use formula
\eqref{62}. $\qquad\Box$

\begin{rem}\label{69}
Notice that for every $\xi=(z,t)\in\h^n\setminus\{0\}$ one has
\begin{eqnarray}
\nonumber&&\frac{1}{\Hn{\xi}^{Q-2}}\sum_{h=1}^n\Hgrad^2(\check{x}_h)X_hu\big(\check\varphi(\xi)\big)+\Hgrad^2(\check{y}_h)Y_hu\big(\check\varphi(\xi)\big)\\
\nonumber&&\hspace{2cm}=\frac{1}{\Hn{\xi}^{Q+6}}G\Hgrad{u}\otimes\big(2(t^2-|z|^4)Jz+4t|z|^2z\big)\\
\nonumber&&\hspace{2,2cm}+\frac{1}{\Hn{\xi}^{Q+6}}\big(2(t^2-|z|^4)Jz+4t|z|^2z\big)\otimes
 G\Hgrad{u}\\
\nonumber&&\hspace{2,2cm}+\frac{1}{\Hn{\xi}^{Q+6}}GJ\Hgrad{u}\otimes\big(2(|z|^4-t^2)z+4t|z|^2Jz\big)\\
\nonumber&&\hspace{2,2cm}+\frac{1}{\Hn{\xi}^{Q+6}}\big(2(|z|^4-t^2)z+4t|z|^2Jz\big)\otimes GJ\Hgrad{u}\\
\nonumber&&\hspace{2,2cm}+\frac{8}{\Hn{\xi}^{Q+6}}\big(\crochet{GJz}{\Hgrad
u}_{\R^{2n}}z-\crochet{Gz}{\Hgrad
u}_{\R^{2n}}Jz\big)\otimes(|z|^2z-tJz)\\
\nonumber&&\hspace{2,2cm}+\frac{8}{\Hn{\xi}^{Q+6}}\big(\crochet{Gz}{\Hgrad u}_{\R^{2n}}z+\crochet{GJz}{\Hgrad u}_{\R^{2n}}Jz\big)\otimes(tz+|z|^2Jz)\\
\nonumber&&\hspace{2,2cm}-\frac{16(|z|^4-t^2)}{\Hn{\xi}^{Q+10}}\big(\crochet{GJz}{\Hgrad
u}_{\R^{2n}}z-\crochet{Gz}{\Hgrad u}_{\R^{2n}}Jz\big)
 \otimes(|z|^2z+tJz)\\
\nonumber&&\hspace{2,2cm}-\frac{32t|z|^2}{\Hn{\xi}^{Q+10}}\big(\crochet{Gz}{\Hgrad
 u}_{\R^{2n}}z+\crochet{GJz}{\Hgrad
 u}_{\R^{2n}}Jz\big)\otimes(|z|^2z+tJz)
\end{eqnarray}
\begin{eqnarray}
\nonumber&&\hspace{2,2cm}+\frac{1}{\Hn{\xi}^{Q+6}}\big(2(|z|^4-t^2)\crochet{GJz}{\Hgrad u}_{\R^{2n}}+4t|z|^2\crochet{Gz}{\Hgrad u}_{\R^{2n}}\big)I_{2n}\\
\nonumber&&\hspace{2,2cm}+\frac{1}{\Hn{\xi}^{Q+6}}\big(2(t^2-|z|^4)\crochet{Gz}{\Hgrad
u}_{\R^{2n}}+4t|z|^2\crochet{GJz}{\Hgrad u}_{\R^{2n}}\big)J,
\end{eqnarray}
where $G,J\in\text{Mat}(2n,\R)$ are the matrices defined in
\eqref{45} and where it's understood that in the previous equality
$\Hgrad{u}$ is to be evaluated at the point $\check\varphi(\xi)$.
\end{rem}

\begin{cor}
Let $u\in\mathcal{C}^2(\h^n)$, then
\begin{equation}\label{70}
\Delta_{\!_H}u_{\check\varphi}(\xi)=\frac{1}{\Hn{\xi}^{Q+2}}\Delta_{\!_H}u(\check{x},\check{y},\check{t}).
\end{equation}
\end{cor}
\textbf{Proof:} In order to conclude one only needs to compute the
trace of the matrix $\Hgrad^2u_{\check\varphi}$, given in formula
\eqref{67}. $\qquad\Box$

On the formula which is the analogue to \eqref{70} when the CR
inversion $\varphi$ is involved, see also \cite{JerLee1},
\cite{BirPra2} and references therein .

\begin{rem}\label{12}
By relations \eqref{Lapl_tr} and \eqref{70} one has
$$u_\psi^{-\frac{Q+2}{Q-2}}\Delta_{\!_H}u_\psi=\big(u^{-\frac{Q+2}{Q-2}}\Delta_{\!_H}u\big)\circ\psi\qquad\text{in }\h^n$$
for every CR map $\psi:\h^n\rightarrow\h^n$, i.e.
$$F(\cdot,u,\Hgrad{u},\Hgrad^2{u}):=u^{-\frac{Q+2}{Q-2}}\Delta_{\!_H}u$$ satisfies definition \ref{def_conf_inv} and it is thus a CR invariant operator on $\h^n$. In particular, if
$u\in\mathcal{C}^2$ is a positive solution of
\begin{equation}\label{37}
-\Delta_{\!_H}u=u^\frac{Q+2}{Q-2},
\end{equation}
so is $u_\psi$ for any CR map $\psi$ on $\h^n$. Equation \eqref{37}
is related to the CR Yamabe problem on $\h^n$, see e.g.
\cite{JerLee1}.
\end{rem}

\section{\textbf{Proof of Theorem \ref{1}}}\label{47}

\begin{lem}\label{3} Given $\xi_0\in\h^n$, $s\in\R^+$, $V\in\R^{2n}$,
$c\in\R$ and $S\in\mathcal{S}^{2n\times2n}$ a symmetric $2n\times2n$
real matrix, there exists $u\in\mathcal{C}^\infty(\h^n)$ which is
positive and such that
$$u(\xi_0)=s,\qquad\Hgrad{u}(\xi_0)=V,\qquad\Hgrad^2u(\xi_0)=S+cJ.$$
\end{lem}
\textbf{Proof:} Let $w\in\mathcal{C}^\infty(\R^{2n+1})$ be a
positive function such that $$w(0)=s,\qquad \nabla
w(0)=\Big(V,\frac{c}{2}\Big),\qquad\nabla^2w(0)=\left(
                                          \begin{array}{c|c}
                                            S & 0 \\
                                            \hline0 & 1 \\
                                          \end{array}
                                        \right).$$
Now define $u(\xi)=w(\xi_0^{-1}\circ\xi)$. Then by relations
\eqref{grad_tr}, \eqref{Hess_tr} and remark \ref{30} one has
\begin{equation*}
\begin{array}{rcl}
u(\xi_0)\!&\!=\!&\!w(0)=s\\
\Hgrad u(\xi_0)\!&\!=\!&\!\Hgrad w(0)=\nabla_zw(0)=V\\
\Hgrad^2u(\xi_0)\!&\!=\!&\!\Hgrad^2w(0)=\nabla^2_zw(0)+2\frac{\partial
w}{\partial t}(0)J=S+cJ
\end{array}
\end{equation*}
as desired.$\qquad\Box$

\textbf{Proof of Theorem \ref{1}}. We start by showing that a second
order differential operator of the form \eqref{68} defined on $\h^n$
which satisfies the invariance properties $(i)$, $(ii)$ and $(iii)$
of the statement is CR invariant on $\h^n$. Indeed, using formulae
\eqref{grad_tr} and \eqref{Hess_tr}, for every positive function
$u\in C^2(\h^n)$ it's easy to see that on $\h^n$ one has
\begin{eqnarray}
\nonumber
A^{u_{\tau_{\hat\xi}}}(\xi)&=&A^u\big(\tau_{\hat\xi}(\xi)\big)\,\,\,\qquad\qquad\qquad\text{for
every }\hat\xi\in\h^n,\\
\nonumber A^{u_{\rho_{\!_M}}}(\xi)&=&\wt{M}^T\cdot
A^u\big(\rho_{\!_M}(\xi)\big)\cdot\wt{M}\qquad\text{for
every }M\in\mathcal{U}(n),\\
\nonumber
A^{u_{\delta_\lambda}}(\xi)&=&A^u\big(\delta_\lambda(\xi)\big)\,\,\,\qquad\qquad\qquad\text{for
every }\lambda>0,\\
\nonumber A^{u_{\iota}}(\xi)&=&G\cdot A^u\big(\iota(\xi)\big)\cdot
G.
\end{eqnarray}
Moreover, using formulae \eqref{62} and \eqref{67} one can see that
for every $\xi\in\h^n\setminus\{0\}$
$$
A^{u_{\check\varphi}}(\xi)\,\,=\,\,E\cdot
A^{u}\big(\check\varphi(\xi)\big)\cdot E^T,
$$
where $G,\,E\in\text{Mat}(2n,\R)$ are the matrices defined
respectively in \eqref{45} and in \eqref{63}.

The first part of the proof is thus complete. We are now going to
prove that a CR invariant differential operator of the second order
on $\h^n$ necessarily satisfies \eqref{68} and the invariance
properties $(i)$, $(ii)$ and $(iii)$ of the statement.

Let $\xi_0\in\h^n$, $s\in\R^+$, $v\in\R^{2n}$ and
$U\in\mathcal{S}^{2n\times2n}\oplus J\R$. Now consider a positive
function $\phi\in\mathcal{C}^\infty(\h^n)$ such that
$\phi(\xi_0)=s$, $\Hgrad{\phi}(\xi_0)=v$ and
$\Hgrad^2\phi(\xi_0)=U$, see lemma \ref{3}. Now use relation
\eqref{conf_inv} with the CR transformation
$$\psi(\xi)=\tau_{\xi_0}(\xi)=\xi_0\circ\xi,$$ see also definition \eqref{trasl}, and the function
$\phi$. By relations \eqref{4}, \eqref{grad_tr} and \eqref{Hess_tr}
one has
$$F\big(\xi,\phi(\xi_0\circ\xi),\Hgrad\phi(\xi_0\circ\xi),\Hgrad^2\phi(\xi_0\circ\xi)\big)=
F\big(\xi_0\circ\xi,\phi(\xi_0\circ\xi),\Hgrad\phi(\xi_0\circ\xi),\Hgrad^2\phi(\xi_0\circ\xi)\big).$$
Evaluating this equality in $\xi=0$ yields
\begin{equation}\label{10}
\begin{array}{l}
F(0,s,v,U)=F\big(0,\phi(\xi_0),\Hgrad\phi(\xi_0),\Hgrad^2\phi(\xi_0)\big)\\
\hspace{2,1cm}=F\big(\xi_0,\phi(\xi_0),\Hgrad\phi(\xi_0),\Hgrad^2\phi(\xi_0)\big)=F(\xi_0,s,v,U).
\end{array}
\end{equation}
Since this holds for every $\xi_0\in\h^n$, $F(\xi_0,s,v,U)$ does not
depend explicitly on $\xi_0$. From now on we will write $F(s,v,U)$
in place of $F(\xi_0,s,v,U)$.

Now let $\phi\in\mathcal{C}^\infty(\h^n)$ be a positive function
such that $\phi(0)=s$, $\Hgrad\phi(0)=0$ and $\Hgrad^2\phi(0)=U$ and
consider the CR transformation
$$\psi(\xi)=\rho_{\!_M}(\xi)=(Bx-Cy,By+Cx,t)$$ for a matrix
$M=B+iC\in\mathcal{U}(n)$, see also definition \eqref{rot}. Using
relation \eqref{conf_inv}, by formulae \eqref{4}, \eqref{grad_tr}
and \eqref{Hess_tr} we get
\begin{equation*}
\begin{array}{l}
F\big(\phi\big(\rho_{\!_M}(\xi)\big),\,\wt{M}^T\cdot\Hgrad\phi\big(\rho_{\!_M}(\xi)\big),\,
\wt{M}^T\cdot\Hgrad^2\phi\big(\rho_{\!_M}(\xi)\big)\cdot\wt{M}\big)\\
\hspace{7cm}=F\big(\phi(\rho_{\!_M}(\xi)),\Hgrad\phi(\rho_{\!_M}(\xi)),\Hgrad^2\phi(\rho_{\!_M}(\xi))\big).
\end{array}
\end{equation*}
Evaluating such equality in $\xi=0$ yields
\begin{equation}\label{M_tilda}
\begin{array}{rcl}
F\big(s,0,\wt{M}^TU\wt{M}\big)&=&F\big(\phi(0),\,\wt{M}^T\cdot\Hgrad\phi(0),\,
\wt{M}^T\cdot\Hgrad^2\phi(0)\cdot\wt{M}\big)\\
&=&F\big(\phi(0),\Hgrad\phi(0),\Hgrad^2\phi(0)\big)\\
&=&F\big(s,0,U\big).
\end{array}
\end{equation}

Consider again a positive function $\phi\in\mathcal{C}^\infty(\h^n)$
 such that $\phi(0)=s$, $\Hgrad\phi(0)=0$ and $\Hgrad^2\phi(0)=U$
and the CR transformation
$$\iota(\xi)=\iota(x,y,t)=(x,-y,-t),$$
see definition \eqref{invers}. Using relations \eqref{conf_inv},
\eqref{4}, \eqref{grad_tr} and \eqref{Hess_tr} as in the previous
cases and evaluating the resulting equality in $0\in\h^n$ gives
\begin{equation}\label{G}
F\big(s,0,GUG\big)=F\big(s,0,U\big).
\end{equation}

Let $\phi\in\mathcal{C}^\infty(\h^n)$ be a positive function such
that $\phi(0)=s$, $\Hgrad\phi(0)=0$ and $\Hgrad^2\phi(0)=U$, let
$Q=2n+2$ and
\begin{equation*}
\psi(\xi)=\delta_{\displaystyle\!
s^{-\frac{2}{Q-2}}}(\xi)=\Big(s^{-\frac{2}{Q-2}}x,\,s^{-\frac{2}{Q-2}}y,\,s^{-\frac{4}{Q-2}}t\Big),
\end{equation*}
see definition \eqref{dil}. By  \eqref{4}, \eqref{grad_tr} and
\eqref{Hess_tr}, relation \eqref{conf_inv} yields
\begin{equation*}
\begin{array}{l}
F\left(s^{-1}\phi\Big(\delta_{\displaystyle\!s^{-\frac{2}{Q-2}}}(\xi)\Big),\,s^{-\frac{Q}{Q-2}}\Hgrad\phi\Big(
\delta_{\displaystyle\!s^{-\frac{2}{Q-2}}}(\xi)\Big),\,s^{-\frac{Q+2}{Q-2}}\Hgrad^2\phi\Big(\delta_{\displaystyle\!s^{-\frac{2}{Q-2}}}(\xi)\Big)\right)\\
\hspace{4cm}=F\left(\phi\Big(\delta_{\displaystyle\!s^{-\frac{2}{Q-2}}}(\xi)\Big),\Hgrad\phi\Big(\delta_{\displaystyle\!
s^{-\frac{2}{Q-2}}}(\xi)\Big),\Hgrad^2\phi\Big(\delta_{\displaystyle\!s^{-\frac{2}{Q-2}}}(\xi)\Big)\right).
\end{array}
\end{equation*}
Evaluating this equality in $\xi=0$ gives
\begin{equation}\label{8}
F\Big(1,0,s^{-\frac{Q+2}{Q-2}}U\Big)=F\big(s,0,U\big).
\end{equation}

Now let $v=(p,q)\in\R^n\times\R^n$, $v\neq0$ and let
$\xi_0:=(x_0,y_0,t_0)$ with
$$x_0=-\frac{(Q-2)s}{|v|^2}p,\quad y_0=-\frac{(Q-2)s}{|v|^2}q,\quad t_0=0$$
and define
$$\lambda:=\Hn{\xi_0}=\big(|x_0|^2+|y_0|^2\big)^\frac{1}{2}=\frac{(Q-2)s}{|v|}.$$
Then consider a positive function $\phi\in\mathcal{C}^\infty(\h^n)$
such that $\phi(\xi_0)=s$, $\Hgrad\phi(\xi_0)=v$ and
$\Hgrad^2\phi(\xi_0)=U$ and the CR transformation
\begin{equation}\label{46}
\psi(\xi)=\varphi\circ\iota\circ\delta_{\lambda^{-2}}(\xi)=
\big(\lambda^2\check{x},\,\lambda^2\check{y},\,\lambda^4\check{t}\big),
\end{equation}
with
\begin{equation}\label{CR-bis}
\check{x}=-\frac{xt+y|z|^2}{\Hn{\xi}^4},\quad\check{y}=\frac{yt-x|z|^2}{\Hn{\xi}^4},\quad\check{t}=\frac{t}{\Hn{\xi}^4},
\end{equation}
see definitions \eqref{dil}, \eqref{invers} and \eqref{CR}. Then
one has $\phi_\psi(\xi)=\frac{\lambda^{Q-2}}{\Hn{\xi}^{Q-2}}
\phi\big(\lambda^2\check{x},\,\lambda^2\check{y},\,\lambda^4\check{t}\big)$,
and by relation \eqref{conf_inv} evaluated in $\psi^{-1}(\xi_0)$ one
has
\begin{eqnarray}
\nonumber F\big(\phi_\psi(\psi^{-1}(\xi_0)),\Hgrad
\phi_\psi(\psi^{-1}(\xi_0)),\Hgrad^2\phi_\psi(\psi^{-1}(\xi_0))\big)\!&\!=\!&\!
F\big(\phi(\xi_0),\Hgrad\phi(\xi_0),\Hgrad^2\phi(\xi_0)\big)\\
\nonumber\!&\!=\!&\!F(s,v,U).
\end{eqnarray}
Now by Lemma \ref{52} in the appendix we have that
$\phi_\psi(\psi^{-1}(\xi_0))=s$,
$\Hgrad\phi_\psi(\psi^{-1}(\xi_0))=0$ and that
\begin{eqnarray}
\nonumber\Hgrad^2\phi_\psi(\psi^{-1}(\xi_0))\!\!&\!\!=\!\!&\!\!G\bigg[-\frac{Q}{Q-2}s^{-1}Jv\otimes
    Jv+\frac{2}{Q-2}s^{-1}v\otimes v+\frac{1}{Q-2}s^{-1}|v|^2I_{2n}\\
\nonumber\!&\!\!&\!\hspace{0,5cm}+J^TUJ+\frac{4}{|v|^4}\Big(\crochet{v}{Uv}_{\R^{2n}}Jv\otimes
    Jv+\crochet{Jv}{UJv}_{\R^{2n}}v\otimes v\\
\label{9}\!&\!\!&\!\hspace{0,5cm}-\crochet{Jv}{\!Uv}_{\R^{2n}}v\otimes
    Jv-\crochet{v}{\!UJv}_{\R^{2n}}\!Jv\otimes v\Big)+\frac{2}{|v|^2}\Big(Jv\otimes
    J^TU^Tv\\
\nonumber\!&\!\!&\!\hspace{0,5cm}-v\otimes J^TU^TJv+J^TUv\otimes
   Jv-J^TUJv\otimes v\Big)\bigg]G,
\end{eqnarray}
with $G$, $J$ being defined as in \eqref{45} and where we recall
that $\crochet{v_1}{v_2}_{\R^{2n}}$ denotes the scalar product of
the vectors $v_1,v_2\in\R^{2n}$. Then
\begin{eqnarray}
\label{6}F\big(s,0,\Hgrad^2\phi_\psi(\psi^{-1}(\xi_0))\big)\!&\!=\!&\!F(s,v,U).
\end{eqnarray}

Now consider the matrix $E\in\mathcal{O}(2n)$ defined in \eqref{63}.
It's not difficult to see that, when evaluated at the point
$\psi^{-1}(\xi_0)$,
\begin{eqnarray}
\nonumber E\!\!&\!\!=\!\!&\!\!G\left(\frac{2}{|v|^2}
   \big(Jv\otimes v-v\otimes Jv\big)+J^T\right).
\end{eqnarray}
Notice that we can apply results \eqref{M_tilda} and \eqref{G} to
relation \eqref{6} using the orthogonal matrix $E$, since it can be
written in the form $G\wt{M}$ with
$$M=B+iC\,\,\in\,\mathcal{U}(n), \qquad \wt{M}=\left(
                 \begin{array}{cc}
            B&-C\\
            C&B
                 \end{array}
               \right)$$
and $$B=\frac{2}{|v|^2}(q\otimes p-p\otimes q),\qquad
C=-\frac{2}{|v|^2}(p\otimes p+q\otimes q)+I_n.$$ We conclude that
\begin{equation}\label{5}
F\big(s,v,U\big)=F\big(s,0,\Hgrad^2\phi_\psi(\psi^{-1}(\xi_0))\big)=F\big(s,0,E^T\,\Hgrad^2\phi_\psi(\psi^{-1}(\xi_0))\,E\big).
\end{equation}
An easy calculation then yields
\begin{eqnarray}
\nonumber
E^T\,\Hgrad^2\phi_\psi(\psi^{-1}(\xi_0))\,E\!\!&\!\!=\!\!&\!\!-\frac{Q}{Q-2}s^{-1}v\otimes
    v+\frac{2}{Q-2}s^{-1}Jv\otimes Jv+\frac{1}{Q-2}s^{-1}|v|^2I_{2n}+U
\end{eqnarray}
and hence from relation \eqref{5} we get
\begin{equation*}
\begin{array}{rcl}
F\big(s,v,U\big)\!\!&\!\!=\!\!&\!\!F\left(s,0,-\frac{Q}{Q-2}s^{-1}v\otimes
    v+\frac{2}{Q-2}s^{-1}Jv\otimes Jv+\frac{1}{Q-2}s^{-1}|v|^2I_{2n}+U\right).
\end{array}
\end{equation*}
This equality trivially holds also in the case $v=0$, then for every
$s\in\R^+$, $v\in\R^{2n}$, $U\in\mathcal{S}^{2n\times2n}\oplus J\R$ by
\eqref{8} one has
\begin{equation}\label{11}
\begin{array}{l}
\hspace{1.3cm}F\big(s,v,U\big)=\\
\hspace{1.5cm}F\left(\!1,0,-\frac{Q}{Q-2}s^{-\frac{2Q}{Q-2}}v\otimes
    v\!+\!\frac{2}{Q-2}s^{-\frac{2Q}{Q-2}}Jv\otimes
    Jv\!+\!\frac{1}{Q-2}s^{-\frac{2Q}{Q-2}}|v|^2I_{2n}\!+\!s^{-\frac{Q+2}{Q-2}}U\!\right).
\end{array}
\end{equation}

Now let $\lambda>0$ and consider the CR map $\psi$ defined in
\eqref{46}. Let $\xi_0=(0,0,\lambda^2)\in\h^n$ and pick a positive
function $\phi\in\mathcal{C}^\infty(\h^n)$ such that
$\phi(\xi_0)=1$, $\Hgrad\phi(\xi_0)=0$ and $\Hgrad^2\phi(\xi_0)=U$.

By Lemma \ref{72} in the appendix we have that
$\phi_\psi(\psi^{-1}(\xi_0))=1$,
$\Hgrad\phi_\psi(\psi^{-1}(\xi_0))=0$ and
$\Hgrad^2\phi_\psi(\psi^{-1}(\xi_0))=-\frac{Q-2}{\lambda^2}J+GUG.$
Since $F$ is CR invariant on $\h^n$, we have
\begin{eqnarray}
\nonumber
F(1,0,U)\!&\!=\!&\!F\big(\phi(\xi_0),\Hgrad\phi(\xi_0),\Hgrad^2\phi(\xi_0)\big)\\
\nonumber\!&\!=\!&\!F\big(\phi_\psi(\psi^{-1}(\xi_0)),\Hgrad
\phi_\psi(\psi^{-1}(\xi_0)),\Hgrad^2\phi_\psi(\psi^{-1}(\xi_0))\big)\\
\nonumber\!&\!=\!&\!F\Big(1,0,-\frac{Q-2}{\lambda^2}J+GUG\Big).
\end{eqnarray}
We now use \eqref{G} together with the relation
$$G\cdot\Big[-\frac{Q-2}{\lambda^2}J+GUG\Big]\cdot G\,\,=\,\,\frac{Q-2}{\lambda^2}J+U$$ to
conclude that for every $\lambda>0$ one has
$$F(1,0,U)=F\Big(1,0,\frac{Q-2}{\lambda^2}J+U\Big).$$
Since $\lambda>0$ is arbitrary, we have that for every $\alpha>0$
\begin{equation}\label{76}
F(1,0,U)=F\Big(1,0,\alpha J+U\Big).
\end{equation}

If we consider $\lambda>0$, the CR map $\psi$ defined in \eqref{46},
the point $\xi_0=(0,0,-\lambda^2)\in\h^n$ and a positive function
$\phi\in\mathcal{C}^\infty(\h^n)$ such that $\phi(\xi_0)=1$,
$\Hgrad\phi(\xi_0)=0$ and $\Hgrad^2\phi(\xi_0)=U$, we can use Lemma
\ref{73} and repeat the above argument. Thus we obtain
\begin{equation}\label{75}
F(1,0,U)=F\Big(1,0,-\alpha J+U\Big)
\end{equation}
for every $\alpha>0$. From equations \eqref{76} and \eqref{75} then
we get
\begin{equation}\label{77}
F(1,0,U)=F\Big(1,0,\alpha J+U\Big)
\end{equation}
for every $\alpha\in\R$. Now notice that
$$U=\Big[\frac{U+U^T}{2}\Big]+\Big[\frac{U-U^T}{2}\Big]=\Big[\frac{U+U^T}{2}\Big]+\alpha J$$ for a suitable
$\alpha\in\R$. Then by relations \eqref{11} and \eqref{77} we
finally get that for every $s\in\R^+$, $v\in\R^{2n}$,
$U\in\mathcal{S}^{2n\times2n}\oplus J\R$ one has
\begin{equation*}
F\big(s,v,U\big)=F\left(1,0,-\frac{Q-2}{2}A\big(s,v,U\big)\right),
\end{equation*}
with
\begin{eqnarray}
\nonumber
   A\big(s,v,U\big)&:=&\frac{2Q}{(Q-2)^2}s^{-\frac{2Q}{Q-2}}v\otimes
   v-\frac{4}{(Q-2)^2}s^{-\frac{2Q}{Q-2}}Jv\otimes Jv\\
\nonumber&&\hspace{2cm}-\frac{2}{(Q-2)^2}s^{-\frac{2Q}{Q-2}}|v|^2I_{2n}-\frac{2}{(Q-2)}s^{-\frac{Q+2}{Q-2}}\Big[\frac{U+U^T}{2}\Big].
\end{eqnarray}
Theorem \ref{1} is now proved.$\qquad\Box$ \vspace{0,1cm}

\begin{rem}
Notice that for every positive function
$u\in\mathcal{C}^2(\h^n)$ one has
$$\text{trace}\big(A^u\big)=-\frac{2}{Q-2}u^{-\frac{Q+2}{Q-2}}\Delta_{\!_H}u$$
which, modulo the harmless constant $-\frac{2}{Q-2}$, is the example
recalled in remark \ref{12}.
\end{rem}

\section{\textbf{Proof of Theorem \ref{13}}}\label{56}

We start this section with a lemma which will be needed in the
course of the proof of Theorem \ref{13}.

For $u\in\mathcal{C}^2(\Omega)$, $u>0$ define
$\phi:=u^{-\frac{2}{Q-2}}$. Then $\phi\in\mathcal{C}^2(\Omega)$,
$\phi>0$ and one has
\begin{equation}\label{16}
\begin{array}{l}
A^u=A_\phi:=\phi\Hhesss\phi-\frac{1}{2}\big|\Hgrad\phi\big|^2I_{2n}-J\Hgrad\phi\otimes
J\Hgrad\phi.
\end{array}
\end{equation}

\begin{lem}\label{17}
Let $\Omega\subset\h^n$ be a bounded open domain,
$\phi\in\mathcal{C}^2(\Omega)$ and $\phi>0$ in $\Omega$. Let
$\eta(\xi)=\eta(z,t):=e^{\delta|z|^2}$ for $\delta>0$. Then there
exists $\bar{\delta}>0$, depending only on $\sup_\Omega|z|$, such
that $\forall\,\delta\in(0,\bar\delta)$ and $\forall\,\e>0$ one has
\begin{equation}\label{23}
A_{\phi+\e\eta}\geq\Big(1+\e\frac{\eta}{\phi}\Big)A_\phi+\e\delta\eta\phi
I_{2n}\qquad\text{in }\Omega.
\end{equation}
\end{lem}

\textbf{Proof:} For every $\phi,\eta\in\mathcal{C}^2(\Omega)$ with
$\phi,\eta>0$ in $\Omega$ and every $\e>0$ one has
\begin{eqnarray}
\nonumber
A_{\phi+\e\eta}\!&\!=\!&\!(\phi+\e\eta)\Hhesss(\phi+\e\eta)-\frac{1}{2}\big|\Hgrad(\phi+\e\eta)\big|^2I_{2n}-
    J\Hgrad(\phi+\e\eta)\otimes J\Hgrad(\phi+\e\eta)\\
\nonumber\!&\!=\!&\!A_\phi+\e\Big(\phi\Hhesss\eta+\eta\Hhesss\phi-\crochet{\Hgrad\phi}{\Hgrad\eta}_{\R^{2n}}I_{2n}\\
\nonumber\!&\!\!&\hspace{4cm}-J\Hgrad\phi\otimes
    J\Hgrad\eta-J\Hgrad\eta\otimes J\Hgrad\phi\Big)+\e^2A_\eta.
\end{eqnarray}
By definition \eqref{16} then one gets
\begin{eqnarray}
\label{19}
\,\,\,\,\,A_{\phi+\e\eta}\!&\!=\!&\!\Big(1+\e\frac{\eta}{\phi}\Big)A_\phi+\e\bigg(\phi\Hhesss\eta+
    \frac{\big|\Hgrad\phi\big|^2}{2\phi}\eta I_{2n}
    +\frac{\eta}{\phi}J\Hgrad\phi\otimes J\Hgrad\phi\\
\nonumber\!&\!\!&\,\,-\crochet{\Hgrad\phi}{\Hgrad\eta}_{\R^{2n}}I_{2n}-J\Hgrad\phi\otimes
    J\Hgrad\eta-J\Hgrad\eta\otimes J\Hgrad\phi\bigg)+\e^2A_\eta.
\end{eqnarray}
Now notice that
\begin{equation*}
\Hgrad\!\Big(\frac{\eta}{\phi}\Big)\!\otimes\!\Hgrad\!\Big(\frac{\eta}{\phi}\Big)=\frac{1}{\phi^2}
    \Hgrad\eta\otimes\Hgrad\eta+\frac{\eta}{\phi^3}\Big(\frac{\eta}{\phi}\Hgrad\phi\otimes\Hgrad\phi-
    \Hgrad\phi\otimes\Hgrad\eta-\Hgrad\eta\otimes\Hgrad\phi\Big),
\end{equation*}
then it follows that
\begin{eqnarray}
\label{18}&&\frac{\phi^3}{\eta}\left(J\Hgrad\!\Big(\frac{\eta}{\phi}\Big)\!\otimes\!J\Hgrad\!\Big(\frac{\eta}{\phi}\Big)-
    \frac{1}{\phi^2}J\Hgrad\eta\otimes J\Hgrad\eta\right)=\\
\nonumber&&\hspace{2,5cm}\left(\frac{\eta}{\phi}J\Hgrad\phi\otimes
    J\Hgrad\phi-J\Hgrad\phi\otimes J\Hgrad\eta-J\Hgrad\eta\otimes J\Hgrad\phi\right).
\end{eqnarray}
Inserting relation \eqref{18} into \eqref{19} we get
\begin{eqnarray}
\label{20}
A_{\phi+\e\eta}\!&\!=\!&\!\Big(1+\e\frac{\eta}{\phi}\Big)A_\phi+\e\bigg(\phi\Hhesss\eta+
    \frac{\big|\Hgrad\phi\big|^2}{2\phi}\eta I_{2n}-\crochet{\Hgrad\phi}{\Hgrad\eta}_{\R^{2n}}I_{2n}\\
\nonumber\!&\!\!&\!\hspace{2cm}\frac{\phi^3}{\eta}J\Hgrad\!\Big(\frac{\eta}{\phi}\Big)\!\otimes\!J\Hgrad\!
    \Big(\frac{\eta}{\phi}\Big)-\frac{\phi}{\eta}J\Hgrad\eta\otimes J\Hgrad\eta\bigg)+\e^2A_\eta.
\end{eqnarray}
Now let $\eta(\xi)=\eta(z,t):=e^{\delta|z|^2}$, with $\delta>0$ to
be chosen later. Then for every $\xi\in\h^n$ one has
$\eta(\xi)\geq1$ and
\begin{equation}\label{21}
\Hgrad\eta(\xi)=2\delta\eta(\xi)z,\qquad\Hhesss\eta(\xi)=2\delta\eta(\xi)I_{2n}+4\delta^2\eta(\xi)z\otimes
z.
\end{equation}
By relations \eqref{21} and since $$0_{2n}\leq v\otimes
v\leq|v|^2I_{2n}\qquad\text{ for all }v\in\R^{2n},$$ it follows that
if $0<\delta\leq\frac{1}{8}\big(\sup_\Omega|z|\big)^{-2}$ one has
\begin{eqnarray}
\nonumber
    A_\eta\!&\!=\!&\!2\delta\eta^2I_{2n}+4\delta^2\eta^2z\otimes
    z-2\delta^2\eta^2|z|^2I_{2n}-4\delta^2\eta^2Jz\otimes Jz\\
\label{38}
    \!&\!\geq\!&\!2\delta\eta^2I_{2n}-2\delta^2\eta^2|z|^2I_{2n}-4\delta^2\eta^2|z|^2I_{2n}\\
\nonumber\!&\!\geq\!&\!\frac{5}{4}\delta\eta^2I_{2n}
\end{eqnarray}
and hence $A_\eta$ is nonnegative definite in $\Omega$. Moreover one
also has
\begin{equation*}
\begin{array}{l}
\phi\Hhesss\eta+\!\frac{\big|\Hgrad\phi\big|^2}{2\phi}\eta
    I_{2n}-\crochet{\Hgrad\phi}{\!\Hgrad\eta}_{\R^{2n}}I_{2n}+\frac{\phi^3}{\eta}
    J\Hgrad\!\Big(\frac{\eta}{\phi}\Big)\!\otimes\!J\Hgrad\!\Big(\frac{\eta}{\phi}\Big)
    \!-\frac{\phi}{\eta}J\Hgrad\eta\otimes J\Hgrad\eta\\
\hspace{1.1cm}\geq\phi\Hhesss\eta+\frac{\big|\Hgrad\phi\big|^2}{2\phi}\eta
    I_{2n}-\crochet{\Hgrad\phi}{\Hgrad\eta}_{\R^{2n}}I_{2n}
    -\frac{\phi}{\eta}J\Hgrad\eta\otimes J\Hgrad\eta\\
\hspace{1.1cm}=2\delta\eta\phi
    I_{2n}+4\delta^2\eta\phi\,z\otimes z+\frac{\big|\Hgrad\phi\big|^2}{2\phi}\eta I_{2n}
    -2\delta\eta\crochet{\Hgrad\phi}{z}_{\R^{2n}}I_{2n}-4\delta^2\eta\phi Jz\otimes Jz\\
\hspace{1.1cm}\geq2\delta\eta\phi
    I_{2n}+\frac{\big|\Hgrad\phi\big|^2}{2\phi}\eta I_{2n}-\frac{\big|\Hgrad\phi\big|^2}{4\phi}\eta I_{2n}
    -4\delta^2\eta\phi|z|^2I_{2n}-4\delta^2\eta\phi|z|^2I_{2n}\\
\hspace{1.1cm}\geq\bigg(\delta\eta\phi+\frac{\big|\Hgrad\phi\big|^2}{4\phi}\eta\bigg)I_{2n}
\end{array}
\end{equation*}
and hence
\begin{equation}\label{39}
\begin{array}{l}
\phi\Hhesss\eta+\frac{\big|\Hgrad\phi\big|^2}{2\phi}\eta
    I_{2n}-\crochet{\Hgrad\phi}{\!\Hgrad\eta}_{\R^{2n}}I_{2n}\\
\hspace{3cm}+\frac{\phi^3}{\eta}
    J\Hgrad\!\Big(\frac{\eta}{\phi}\Big)\otimes J\Hgrad\!\Big(\frac{\eta}{\phi}\Big)
    -\frac{\phi}{\eta}J\Hgrad\eta\otimes J\Hgrad\eta\geq\delta\eta\phi I_{2n}.
\end{array}
\end{equation}
Then from \eqref{20} by relations \eqref{38} and \eqref{39} we have
$$A_{\phi+\e\eta}\geq\Big(1+\e\frac{\eta}{\phi}\Big)A_\phi+\e\delta\eta\phi I_{2n}\qquad\text{in }\Omega.$$
The proof of the lemma is now complete.$\qquad\Box$

\begin{rem}
Inequality \eqref{23} is equivalent to
\begin{equation*}
\begin{array}{rcl}
A^{\!^{\big[\big(u^{-\frac{2}{Q-2}}+\e\eta\big)^{-\frac{Q-2}{2}}\big]}}&\geq&\Big(1+\e\eta
     u^\frac{2}{Q-2}\Big)A^u+\e\delta\eta u^{-\frac{2}{Q-2}}I_{2n}\qquad\text{in
     }\Omega.
\end{array}
\end{equation*}
\end{rem}
\hspace{0,1cm} We are now ready to start with the proof of Theorem
\ref{13}.

\textbf{Proof of Theorem \ref{13}:} We start by proving part (i) of
the statement.

Let $\phi:=u^{-\frac{2}{Q-2}}$ and $\theta:=w^{-\frac{2}{Q-2}}$.
Then $u(\xi)\geq w(\xi)$ if and only if $\phi(\xi)\leq\theta(\xi)$,
for any $\xi\in\overline{\Omega}$. Hence in particular we have
$\phi\leq\theta$ on $\partial\Omega$. Moreover we also have
\begin{equation}\label{40}
A_\phi=A^u\in\overline{\Sigma},\,\,\,\,\,\,\,\,\,A_\theta=A^w\in\Sigma^c\quad
\text{in }\Omega.
\end{equation}
Now by contradiction suppose there exists $\xi_0\in\Omega$ such that
$u(\xi_0)<w(\xi_0)$, i.e. such that $\phi(\xi_0)>\theta(\xi_0)$.
Multiply $\phi$ by a constant $\alpha_*\in(0,1)$ so that
\begin{equation}\label{26}
\begin{array}{rcl}
\theta\!&\!>\!&\!\alpha_*\phi\,\,\,\,\qquad\text{ on }\partial\Omega,\\
\theta\!&\!\geq\!&\!\alpha_*\phi\,\,\,\,\qquad\text{ in }\Omega,\\
\theta(\xi_1)\!&\!=\!&\!\alpha_*\phi(\xi_1)\quad\text{ for some
}\xi_1\in\Omega.
\end{array}
\end{equation}
One can easily prove that $\alpha_*=\sup\big\{\alpha\in(0,1):\beta
\phi\leq \theta\text{ in }\Omega,\,\forall\,\beta\in(0,\alpha)\big\}$.
 Now use lemma \ref{17}, and let $\eta(\xi)=e^{\delta|z|^2}$ for
 some $\delta>0$ small enough so that for $\e>0$ one has
 \begin{equation}\label{36}
A_{\alpha_*\phi+\e\eta}\geq\Big(1+\e\frac{\eta}{\alpha_*\phi}\Big)A_{\alpha_*\phi}+\e\delta\eta
(\alpha_*\phi)I_{2n}\qquad\text{in }\Omega.
\end{equation}
Choose $\e>0$ small
enough, so that
$$\theta>\alpha_*\phi+\e\eta\qquad\text{on }\partial\Omega.$$ For instance this
can be achieved by choosing
$\displaystyle0<\e<\frac{1}{2}\Big(\inf_{\partial\Omega}(\theta-\alpha_*\phi)\Big)\Big(\sup_{\partial\Omega}\eta\Big)^{-1}$.
Then we have $\theta(\xi_1)<\alpha_*\phi(\xi_1)+\e\eta(\xi_1)$.

Now, in a similar way as we already did with $\alpha_*$ in relations
\eqref{26}, let $\gamma\in(0,1)$ be such that
\begin{equation}\label{27}
\begin{array}{rcl}
\theta\!&\!>\!&\!\gamma\big(\alpha_*\phi+\e\eta\big)\qquad\,\,\qquad\text{ on }\partial\Omega,\\
\theta\!&\!\geq\!&\!\gamma\big(\alpha_*\phi+\e\eta\big)\qquad\,\,\qquad\text{ in }\Omega,\\
\theta(\xi_2)\!&\!=\!&\!\gamma\big(\alpha_*\phi(\xi_2)+\e\eta(\xi_2)\big)\quad\text{
for some }\xi_2\in\Omega,
\end{array}
\end{equation}
where
$\gamma=\sup\big\{c\in(0,1):\beta\big(\alpha_*\phi+\e\eta\big)\leq\theta\text{
in }\Omega,\,\forall\,\beta\in(0,c)\big\}$. Consider the CR map
$\tau_{\xi_2}(\xi)$ and the transformed functions
$\phi_{\tau_{\xi_2}}$, $\theta_{\tau_{\xi_2}}$ and
$\eta_{\tau_{\xi_2}}$, see also definition \eqref{trasl} and
relation \eqref{4}. By relations \eqref{27} then we have
\begin{eqnarray}
\nonumber\theta_{\tau_{\xi_2}}\!&\!\geq\!&\!\gamma\big(\alpha_*\phi_{\tau_{\xi_2}}+\e\eta_{\tau_{\xi_2}}\big)\qquad
     \qquad\text{ in }\tau_{\xi_2}^{-1}\big(\overline{\Omega}\big),\\
\label{35}\theta_{\tau_{\xi_2}}(0)\!&\!=\!&\!\gamma\big(\alpha_*\phi_{\tau_{\xi_2}}(0)+\e\eta_{\tau_{\xi_2}}(0)\big).
\end{eqnarray}
Then we have
\begin{eqnarray}
\label{28}\nabla\theta_{\tau_{\xi_2}}(0)\!&\!=\!&\!\gamma\big(\alpha_*\nabla\phi_{\tau_{\xi_2}}(0)+\e\nabla
    \eta_{\tau_{\xi_2}}(0)\big),\\
\label{29}\nabla^2\theta_{\tau_{\xi_2}}(0)\!&\!\geq\!&\!\gamma\big(\alpha_*\nabla^2\phi_{\tau_{\xi_2}}(0)+\e\nabla^2
    \eta_{\tau_{\xi_2}}(0)\big),
\end{eqnarray}
where we recall that $\nabla$ and $\nabla^2$ denote respectively the
gradient and the Hessian matrix of a regular function in $\R^{2n+1}$. Now recall that by
remark \ref{30} for any function $f\in\mathcal{C}^1(\Omega)$ one has $\Hgrad f(0)=\nabla_zf(0)$, hence by relation
\eqref{28} we have
\begin{equation}\label{31}
\Hgrad\theta_{\tau_{\xi_2}}(0)=\gamma\big(\alpha_*\Hgrad\phi_{\tau_{\xi_2}}(0)+\e\Hgrad\eta_{\tau_{\xi_2}}(0)\big).
\end{equation}
We notice also that \eqref{29} implies
\begin{equation}\label{34}
\Hhesss\theta_{\tau_{\xi_2}}(0)\geq\gamma\big(\alpha_*\Hhesss\phi_{\tau_{\xi_2}}(0)+\e\Hhesss\eta_{\tau_{\xi_2}}(0)\big).
\end{equation}
Indeed for any $z\in\R^{2n}$ let $\zeta:=(z,0)\in\R^{2n+1}$, then we have
\begin{equation*}
\begin{array}{l}
\displaystyle\crochet{z}{\nabla_z^2\theta_{\tau_{\xi_2}}(0)\,z}_{\R^{2n}}=\sum_{i,j=1}^{2n}\big[\nabla_z^2
    \theta_{\tau_{\xi_2}}(0)\big]_{ij}z_iz_j\\
\hspace{5,5cm}\displaystyle=\sum_{i,j=1}^{2n+1}\big[\nabla^2\theta_{\tau_{\xi_2}}(0)\big]_{ij}\zeta_i\zeta_j=
    \crochet{\zeta}{\nabla^2\theta_{\tau_{\xi_2}}(0)\,\zeta}_{\R^{2n+1}}\\
\displaystyle\crochet{z}{\gamma\big(\alpha_*\nabla_z^2\phi_{\tau_{\xi_2}}(0)+\e\nabla_z^2
    \eta_{\tau_{\xi_2}}(0)\big)\,z}_{\R^{2n}}=\crochet{\zeta}{\gamma\big(\alpha_*\nabla^2\phi_{\tau_{\xi_2}}(0)
    +\e\nabla^2\eta_{\tau_{\xi_2}}(0)\big)\,\zeta}_{\R^{2n+1}}
\end{array}
\end{equation*}
Now, by relation \eqref{29}, this in turn implies that
\begin{equation}\label{32}
\nabla_z^2\theta_{\tau_{\xi_2}}(0)\geq\gamma\big(\alpha_*\nabla^2_z\phi_{\tau_{\xi_2}}(0)
    +\e\nabla^2_z\eta_{\tau_{\xi_2}}(0)\big).
\end{equation}
By remark \ref{30} one has $\Hhesss{f}(0)=\nabla_z^2f(0)$ for every
function $f\in\mathcal{C}^2(\Omega)$, so that inequality \eqref{32}
finally implies \eqref{34}.

But then by formulae \eqref{grad_tr} and \eqref{Hess_tr} we have
\begin{eqnarray}
\nonumber\theta(\xi_2)\!\!&\!\!=\!\!&\!\!\theta_{\tau_{\xi_2}}(0)=\gamma\big(\alpha_*\phi_{\tau_{\xi_2}}(0)+\e
    \eta_{\tau_{\xi_2}}(0)\big)=\gamma\big(\alpha_*\phi(\xi_2)+\e\eta(\xi_2)\big)\,\qquad\text{by \eqref{35}},\\
\nonumber\Hgrad\theta(\xi_2)\!\!&\!\!=\!\!&\!\!\Hgrad\theta_{\tau_{\xi_2}}(0)=\gamma\big(\alpha_*\Hgrad
    \phi_{\tau_{\xi_2}}(0)+\e\Hgrad\eta_{\tau_{\xi_2}}(0)\big)\qquad\qquad\qquad\qquad\text{by \eqref{31}}\\
\nonumber\!\!&\!\!\!\!&\!\!\hspace{4.6cm}=\gamma\big(\alpha_*\Hgrad\phi(\xi_2)+\e\Hgrad\eta(\xi_2)\big),\\
\nonumber\Hhesss\theta(\xi_2)\!\!&\!\!=\!\!&\!\!\Hhesss\theta_{\tau_{\xi_2}}(0)\geq\gamma\big(\alpha_*\Hhesss
    \phi_{\tau_{\xi_2}}(0)+\e\Hhesss\eta_{\tau_{\xi_2}}(0)\big)\qquad\qquad\qquad\qquad\,\text{by \eqref{34}}\\
\nonumber\!\!&\!\!\!\!&\!\!\hspace{4.6cm}=\gamma\big(\alpha_*\Hhesss\phi(\xi_2)+\e\Hhesss\eta(\xi_2)\big).
\end{eqnarray}
This in turn implies that at $\xi_2\in\Omega$ one has
\begin{eqnarray}
\nonumber
A_\theta\!\!&\!\!=\!\!&\!\!\theta\Hhesss\theta-\frac{1}{2}\big|\Hgrad\theta\big|^2I_{2n}-
    J\Hgrad\theta\otimes J\Hgrad\theta\\
\nonumber\!\!&\!\!\geq\!\!&\!\!A_{\gamma(\alpha_*\phi+\e\eta)}\\
\label{41}\!\!&\!\!=\!\!&\!\!\gamma^2A_{\alpha_*\phi+\e\eta}\\
\nonumber \!\!&\!\!\geq\!\!&\!\!\gamma^2\Big(\Big(1+\e\frac{\eta}{\alpha_*\phi}\Big)A_{\alpha_*\phi}+\e\delta\eta
     (\alpha_*\phi)I_{2n}\Big)\qquad\text{by \eqref{36}}\\
\nonumber \!\!&\!\!=\!\!&\!\!\gamma^2\alpha_*^2\Big(1+\e\frac{\eta}{\alpha_*\phi}\Big)A_\phi+\e\delta\gamma^2\alpha_*\eta
     \phi I_{2n}.
\end{eqnarray}
Since $A_\phi\in\overline\Sigma$ and since $\gamma^2\alpha_*^2\Big(1+\e\frac{\eta}{\alpha_*\phi}\Big)>0$ in $\Omega$,
by condition \eqref{14} we get
\begin{equation}\label{42}
\gamma^2\alpha_*^2\Big(1+\e\frac{\eta}{\alpha_*\phi}\Big)A_\phi\in\overline\Sigma.
\end{equation}
Then in $\xi_2\in\Omega$ we have $$A_\theta=\Big(A_\theta-\gamma^2\alpha_*^2\Big(1+\e\frac{\eta}{\alpha_*\phi}\Big)A_\phi
\Big)+\gamma^2\alpha_*^2\Big(1+\e\frac{\eta}{\alpha_*\phi}\Big)A_\phi:=B+cA_\phi,$$
with $cA_\phi\in\overline\Sigma$ by \eqref{42} and with $B\in\mathcal{S}^{2n\times2n}$, $B\geq\e\delta\gamma^2\alpha_*\eta\phi I_{2n}>0$ by \eqref{41}.

By condition \eqref{14} it follows that $A_\theta\in\Sigma$ when evaluated in $\xi_2\in\Omega$, which contradicts our
hypothesis $A_\theta\in\Sigma^c$ in $\Omega$, see condition \eqref{40}. Then we have $u\geq w$ in $\overline\Omega$.

Part (i) of the statement of the theorem is thus proved. Now we turn
our attention to part (ii).

Consider again $\phi:=u^{-\frac{2}{Q-2}}$ and $\theta:=w^{-\frac{2}{Q-2}}$. Then $A_\phi\in\overline\Sigma$ and
$A_\theta\in\Sigma^c$. Since $u>w$ on $\partial\Omega$, we have
$\theta>\phi$ on $\partial\Omega$. By part (i) we have $u\geq w$ in $\overline\Omega$, and hence $\theta\geq\phi$ in
$\overline\Omega$.

Now suppose by contradiction that there exists $\xi_1\in\Omega$ such that $u(\xi_1)=w(\xi_1)$, i.e.
$\theta(\xi_1)=\phi(\xi_1)$. Then we have condition \eqref{26}, this time with $\alpha_*=1$.

The proof now proceeds as in the previous case (i), where one has just to substitute $\alpha_*=1$ in all the calculations.
Then we can conclude that at some point $\xi_2\in\Omega$ one has $A_\theta\in\Sigma$, which again contradicts our
hypothesis $A_\theta\in\Sigma^c$ in $\Omega$. Thus $u>w$ in $\overline\Omega$, and the proof of the theorem is
now complete.$\qquad\Box$

\vspace{0,1cm}

\begin{rem}
Notice that, by choosing
$$\Sigma:=\Big\{A\in\mathcal{S}^{2n\times2n}\,\Big|\,\text{tr}A>0\Big\}$$ in Theorem \ref{13}, we have the following
corollary:

\vspace{0,2cm} Let
$u,w\in\mathcal{C}^2(\Omega)\cup\mathcal{C}^0(\overline{\Omega})$
with $\Omega\subset\h^n$ a bounded open domain. Suppose that $u,w>0$
in $\overline{\Omega}$, and that $\Delta_{\!_H}u\leq0$,
$\Delta_{\!_H}w\geq0$ in $\Omega$. Then
\begin{itemize}
\item[i)] if $u\geq w$ on
$\partial\Omega$, $u\geq w$ in $\overline{\Omega}$,
\item[ii)] if $u>w$ on
$\partial\Omega$, $u> w$ in $\overline{\Omega}$.
\end{itemize}
\vspace{0,1cm} This is also a consequence of the weak maximum
principle applied to the sublaplacian on the Heisenberg group.
\end{rem}

\section{\textbf{Proof of Theorem \ref{94}}}\label{95}

We start this section with some results that will be needed in the
course of the proof of Theorem \ref{94}.

\begin{thm}\label{93}
Let $n\geq1$, $Q=2n+2$ and consider $D_1(0)\subset\h^n$,
$u\in\mathcal{C}^2\big(D_1(0)\setminus\{0\}\big)$ such that
$$\Delta_{\!_H}u\leq0\qquad\text{in }D_1(0)\setminus\{0\}.$$
Let $f,g:D_1(0)\rightarrow\R$ be bounded functions such that
\begin{itemize}
\item[i)] $f(0)=g(0)$,
\item[ii)] $f,g$ are differentiable in $0$ and
$\Hgrad{f}(0)\neq\Hgrad{g}(0)$,
\item[iii)] $u(\xi)\geq f(\xi)$, $u(\xi)\geq g(\xi)$ for every $\xi\in
D_1(0)\setminus\{0\}$.
\end{itemize}
Then one has
$\displaystyle\lim_{R\rightarrow0^+}\inf_{D_R(0)}u>f(0)$.
\end{thm}
\textbf{Proof:} For every $\xi\in D_1(0)\setminus\{0\}$ define
\begin{eqnarray}
\nonumber\wt{u}(\xi)&:=&u(\xi)-f(0)-\crochet{\nabla
f(0)}{\xi}_{\R^{2n+1}},
\end{eqnarray}
and for $\xi\in D_1(0)$ define
\begin{eqnarray}
\nonumber\wt{f}(\xi)&:=&f(\xi)-f(0)-\crochet{\nabla
f(0)}{\xi}_{\R^{2n+1}},\\
\nonumber\wt{g}(\xi)&:=&g(\xi)-f(0)-\crochet{\nabla
f(0)}{\xi}_{\R^{2n+1}},
\end{eqnarray}
where $\nabla$ denotes the usual gradient of a function defined in a
domain of $\R^{2n+1}$. Then $\wt{u}$, $\wt{f}$, $\wt{g}$ satisfy all
of the hypotheses of the theorem, moreover with
$$\wt{f}(0)=\wt{g}(0)=0,\,\,\,\,\,\,\nabla\wt{f}(0)=0,\,\,\,\,\,\,\Hgrad\wt{f}(0)=0,\,\,\,\,\,\,\Hgrad\wt{g}(0)\neq0.$$
If we prove that the result of the theorem holds for the functions
$\wt{u}$, $\wt{f}$, $\wt{g}$, that is if we have
$$\lim_{R\rightarrow0}\inf_{D_R(0)}\wt{u}>0=\wt{f}(0),$$ then we also have
$$\lim_{R\rightarrow0}\Big(\inf_{D_R(0)}u-f(0)-\inf_{D_R(0)}\crochet{\nabla
f(0)}{\xi}_{\R^{2n+1}}\Big)
\geq\lim_{R\rightarrow0}\inf_{D_R(0)}\wt{u}>0.$$ Thus we obtain the
desired result on $u$, i.e.
$\displaystyle\lim_{R\rightarrow0}\inf_{D_R(0)}u>f(0).$

Hence without loss of generality we can assume that the functions
$u,\,f,\,g$ also satisfy
\begin{equation}\label{79}
f(0)=g(0)=0,\,\,\,\,\,\,\nabla f(0)=0,\,\,\,\,\,\,\Hgrad
f(0)=0,\,\,\,\,\,\,\Hgrad g(0)\neq0,
\end{equation}
and thus we have to prove that
$\displaystyle\lim_{R\rightarrow0}\inf_{D_R(0)}u>0=f(0).$

Conditions \eqref{79} in particular imply that as
$|\xi|\rightarrow0$ one has
\begin{equation*}
f(\xi)=o(|\xi|),\,\,\,\,\,\,\,\,\,\,\,\,\,g(\xi)=\crochet{\nabla
g(0)}{\xi}_{\R^{2n+1}}+o(|\xi|),
\end{equation*}
where $|\cdot|$ denotes the usual Euclidean norm. Then in $D_1(0)$
we have
\begin{equation}\label{80}
u(\xi)\geq-h(\xi),\,\,\,\,\,\,\,\,\,\,u(\xi)\geq\crochet{\nabla
g(0)}{\xi}_{\R^{2n+1}}-h(\xi):=\crochet{\zeta}{\xi}_{\R^{2n+1}}-h(\xi)
\end{equation}
with $\zeta=(z_1,t_1):=\nabla g(0)\in\R^{2n}\times\R$,
$z_1=\Hgrad{g}(0)\neq0$ by \eqref{79}, and with $h(\xi)\geq0$
bounded and $h(\xi)=o(|\xi|)$ as $|\xi|\rightarrow0$. Now define on
$D_1(0)\setminus\{0\}$ the function
$$u_\lambda(\xi):=\frac{1}{\lambda}u\big(\delta_\lambda(\xi)\big),$$
with $0<\lambda\ll1$ and $\delta_\lambda(\xi)$ defined by relation
\eqref{dil}. Notice that, as $\lambda$ tends to $0$, we have on
$\overline{D_1(0)}$
$$\frac{1}{\lambda}h\big(\delta_\lambda(\xi)\big)=\frac{1}{\lambda}|\delta_\lambda(\xi)|\,o(1)=
\frac{\sqrt{\lambda^2|z|^2+\lambda^4t^2}}{\lambda}\,o(1)\leq|\xi|\,o(1)\leq\sqrt{2}\,o(1).$$
Hence $h\big(\delta_\lambda(\xi)\big)=o(\lambda)$ as
$\lambda\rightarrow0$, uniformly for $\xi\in\overline{D_1(0)}$. By
relations \eqref{80} then we have
\begin{eqnarray}
\label{81}u_\lambda(\xi)\!&\!\geq\!&\!-\frac{1}{\lambda}h\big(\delta_\lambda(\xi)\big)\geq\,-o(1)\qquad\text{on
}\overline{D_1(0)}\setminus\{0\}\\
\label{82}u_\lambda(\xi)\!&\!\geq\!&\!\frac{1}{\lambda}\crochet{\zeta}{\delta_\lambda(\xi)}_{\R^{2n+1}}-\frac{1}{\lambda}
h\big(\delta_\lambda(\xi)\big)\\
\nonumber\!&\!\geq\!&\!\crochet{z_1}{z}_{\R^{2n}}+\lambda
t_1t-o(1)\geq\crochet{z_1}{z}_{\R^{2n}}-0(1)\qquad\text{on
}\overline{D_1(0)}\setminus\{0\},
\end{eqnarray}
uniformly in $\xi=(z,t)$ as $\lambda\rightarrow0$.

By \eqref{81} for every $\e>0$ there exists $\lambda_0>0$ such that
\begin{equation}\label{83}
u_\lambda(\xi)\geq-\e\qquad\text{in
}\overline{D_1(0)}\setminus\{0\}\,\text{ for all }\lambda<\lambda_0.
\end{equation}

Moreover, if we set $\xi_0:=\left(\frac{z_1}{2|z_1|},0\right)$, we
have $\Hn{\xi_0}=|\xi_0|=\frac{1}{2}$ and it's easy to see that
$\overline{D_\frac{1}{4}(\xi_0)}\subset D_1(0)\setminus\{0\}$. By
\eqref{82} on $\overline{D_\frac{1}{4}(\xi_0)}$ then we have
\begin{eqnarray}
\nonumber u_\lambda(\xi)\!&\!\geq\!&\!\crochet{z_1}{z}_{\R^{2n}}-o(1)=\frac{1}{2}|z_1|+\crochet{z_1}{z-\frac{z_1}{2|z_1|}}_{\R^{2n}}-o(1)\\
\label{84}\!&\!\geq\!&\!\frac{1}{2}|z_1|-|z_1|\left|z-\frac{z_1}{2|z_1|}\right|-o(1)\\
\nonumber\!&\!\geq\!&\!\frac{1}{2}|z_1|-|z_1|\Hn{\xi-\xi_0}-o(1)\\
\nonumber\!&\!\geq\!&\!\frac{1}{2}|z_1|-\frac{1}{4}|z_1|-o(1)>c_0,
\end{eqnarray}
with $c_0:=\frac{1}{8}|z_1|>0$, for every positive $\lambda$ which
is smaller than a suitable $\lambda_1>0$. From now on we will assume
that $0<\lambda<\bar{\lambda}:=\min\{\lambda_0,\lambda_1\}$.

Now let $\sigma^{(\e)}$ be the solution of
\begin{equation}\label{85}
\begin{cases}
\Delta_{\!_H}\sigma^{(\e)}=0\qquad\text{ on }
D_1(0)\setminus\overline{D_\frac{1}{4}(\xi_0)}\\
\sigma^{(\e)}=-2\e\qquad\,\,\text{ on }\partial D_1(0)\\
\sigma^{(\e)}=\frac{1}{2}c_0\qquad\,\,\,\,\text{ on }\partial
D_\frac{1}{4}(\xi_0)
\end{cases}
\end{equation}
Since $\Omega:=D_1(0)\setminus\overline{D_\frac{1}{4}(\xi_0)}$ is a
smooth domain and since its boundary is characteristic for
$\Delta_{\!_H}$ only in the north and south poles of the two balls,
i.e. in
$$N_1=(0,1),\,\,\,S_1=(0,-1),\,\,\,N_2=\left(\frac{\zeta_z}{2|\zeta_z|},\frac{1}{16}\right),\,\,\,
S_2=\left(\frac{\zeta_z}{2|\zeta_z|},-\frac{1}{16}\right),$$ problem
\eqref{85} admits a unique solution which is $\mathcal{C}^\infty$ in
the interior of the domain, because the operator is hypoelliptic
since it satisfies H\"{o}rmander's condition (see \cite{Hor}), and
also up to the boundary in any point which is different from
$N_1,S_1,N_2,S_2$ (see \cite{KohNin} and \cite{Jer1}). The unique
solution of problem \eqref{85} is also continuous up to the boundary
in all the points, see \cite{Jer1}.

By the maximum principle, see \cite{Bony}, the solution depends
continuously on the data of the problem and as $\e$ tends to $0$ we
have
\begin{equation}\label{86}
\sup_{\overline{\Omega}}\big|\sigma^{(\e)}-\sigma^{(0)}\big|\rightarrow0,
\end{equation}
where $\sigma^{(0)}$ is the unique solution of
\begin{equation}\label{87}
\begin{cases}
\Delta_{\!_H}\sigma^{(0)}=0\qquad\text{ on }
D_1(0)\setminus\overline{D_\frac{1}{4}(\xi_0)}\\
\sigma^{(0)}=0\qquad\,\,\,\,\,\,\,\,\,\text{ on }\partial D_1(0)\\
\sigma^{(0)}=\frac{1}{2}c_0\qquad\,\,\,\,\text{ on }\partial
D_\frac{1}{4}(\xi_0).
\end{cases}
\end{equation}
By the strong maximum principle, see \cite{Bony}, we have that
$\sigma^{(0)}\geq0$ in $\overline\Omega$ and $\sigma^{(0)}$ cannot
attain its minimum value in $\Omega$. In particular we have
$\sigma^{(0)}(0)>0$. Then by \eqref{86} if we choose $\e>0$ small
enough we have
\begin{equation}\label{88}
\sigma^{(\e)}(0)\geq\frac{1}{2}\sigma^{(0)}(0)>0.
\end{equation}
Notice moreover that by the maximum principle one also has
\begin{equation}\label{89}
-2\e\leq\sigma^{(\e)}\leq\frac{1}{2}c_0\qquad\text{ on
}\overline\Omega.
\end{equation}
Now fix $\e>0$ such that \eqref{88} holds, and for
$0<\lambda<\bar\lambda$ and any $0<r<\frac{1}{8}$ define
\begin{equation}\label{90}
\Theta_\lambda(\xi)=u_\lambda(\xi)+Ar^{Q-2}\big(\Hn\xi^{-(Q-2)}-1\big)-\sigma^{(\e)}(\xi)
\end{equation}
on $D_1(0)\setminus\overline{\big(D_r(0)\cup
D_\frac{1}{4}(\xi_0)\big)}$, with $A>1$ to be fixed later. Notice
that $\overline{D_r(0)}\subset D_1(0)$ and that
$\overline{D_r(0)}\cap\overline{D_\frac{1}{4}(\xi_0)}=\emptyset$.
Moreover one also has $\Delta_{\!_H}\big(\Hn\xi^{-(Q-2)}\big)=0$ in
$\h^n\setminus\{0\}$, since $\Hn\xi^{-(Q-2)}$ is a constant multiple
of the fundamental solution of $\Delta_{\!_H}$ centered at $0$.

Then $\Delta_{\!_H}\Theta_\lambda=\Delta_{\!_H}u_\lambda\leq0$ in
$D_1(0)\setminus\overline{\big(D_r(0)\cup
D_\frac{1}{4}(\xi_0)\big)}$.

Now notice that $\sigma^{(\e)}=-2\e$ on $\partial D_1(0)$, and hence
$\sigma^{(\e)}\leq-\frac{3}{2}\e$ near $\partial D_1(0)$. But then
near $\partial D_1(0)$ by \eqref{83} one has
$$\Theta_\lambda\geq u_\lambda-\sigma^{(\e)}\geq-\e+\frac{3}{2}\e=\frac{1}{2}\e>0.$$

On $\partial D_\frac{1}{4}(\xi_0)$, by \eqref{84} and \eqref{85},
one has
$$\Theta_\lambda\geq u_\lambda-\sigma^{(\e)}\geq c_0-\frac{1}{2}c_0=\frac{1}{2}c_0>0.$$

Finally on $\partial D_r(0)$, by \eqref{90} and \eqref{89}, one has
$$\Theta_\lambda\geq -\e+Ar^{Q-2}\big(r^{-(Q-2)}-1\big)-\frac{1}{2}c_0\geq
A\left(1-\frac{1}{8^{Q-2}}\right)-\e-\frac{1}{2}c_0>0$$ if we choose
$A>1$ large enough. Hence $\Theta_\lambda>0$ on
$\partial\big(D_1(0)\setminus\overline{\big(D_r(0)\cup
D_\frac{1}{4}(\xi_0)\big)}\,\big)$. From the maximum principle then
it follows that $\Theta_\lambda>0$ on the whole domain.

Now fix $\xi\in D_1(0)\setminus\{0\}$ and let
$0<r<\min\big\{\Hn{\xi},\frac{1}{8}\big\}$. For all
$0<\lambda<\bar{\lambda}$ then we have
$$u_\lambda(\xi)+Ar^{Q-2}\big(\Hn\xi^{-(Q-2)}-1\big)-\sigma^{(\e)}(\xi)>0,$$
and by letting $r$ tend to $0$ we get
$u_\lambda(\xi)\geq\sigma^{(\e)}(\xi)$.

Then for all $R\in(0,\lambda)$ one has
$$\inf_{D_R(0)}u(\xi)=\inf_{D_\frac{R}{\lambda}(0)}u\big(\delta_\lambda(\xi)\big)=
\lambda\inf_{D_\frac{R}{\lambda}(0)}u_\lambda(\xi)\geq\lambda\inf_{D_\frac{R}{\lambda}(0)}\sigma^{(\e)}(\xi).$$
Hence by the continuity of $\sigma^{(\e)}$ and by relation
\eqref{88} we get
\begin{eqnarray}
\nonumber\lim_{R\rightarrow0}\inf_{D_R(0)}u(\xi)\geq\lambda
\lim_{R\rightarrow0}\inf_{D_\frac{R}{\lambda}(0)}\sigma^{(\e)}(\xi)=\lambda\,\sigma^{(\e)}(0)
\geq\frac{\lambda}{2}\sigma^{(0)}(0)>0.
\end{eqnarray}
Thus the proof of the theorem is now complete.$\qquad\Box$

Exploiting left translations and dilations, defined respectively in
\eqref{trasl} and \eqref{dil}, one can easily derive the following
Corollary from Theorem \ref{93}.
\begin{cor}\label{91} Let $n\geq1$, $Q=2n+2$, $R>0$, $\xi_0\in\h^n$ and consider $D_R(\xi_0)\subset\h^n$,
$u\in\mathcal{C}^2\big(D_R(\xi_0)\setminus\{\xi_0\}\big)$ such that
$$\Delta_{\!_H}u\leq0\qquad\text{in }D_R(\xi_0)\setminus\{\xi_0\}.$$
Let $f,g:D_R(\xi_0)\rightarrow\R$ be bounded functions such that
\begin{itemize}
\item[i)] $f(\xi_0)=g(\xi_0)$,
\item[ii)] $f,g$ are differentiable in $\xi_0$ and
$\Hgrad{f}(\xi_0)\neq\Hgrad{g}(\xi_0)$,
\item[iii)] $u(\xi)\geq f(\xi)$, $u(\xi)\geq g(\xi)$ for every $\xi\in
D_R(\xi_0)\setminus\{\xi_0\}$.
\end{itemize}
Then one has
$\displaystyle\lim_{r\rightarrow0^+}\inf_{D_r(\xi_0)}u>f(\xi_0)$.
\end{cor}

For any positive function $w\in\mathcal{C}^2\big(D_2(0)\big)$,
$\xi\in D_2(0)$ and $\lambda>0$ let
\begin{equation}\label{96}
w^{\xi,\lambda}(\eta)=w_{\tau_\xi\circ\delta_\lambda}(\eta)=\lambda^\frac{Q-2}{2}w\big(\xi\circ\delta_\lambda(\eta)\big),
\end{equation}
for $\eta\in\h^n$ such that $\xi\circ\delta_\lambda(\eta)\in
D_2(0)$.

\begin{lem}\label{15}
Assume $u\in\mathcal{C}^2\big(D_2(0)\setminus\{0\}\big)$,
$w\in\mathcal{C}^2\big(D_2(0)\big)$, $u>w$ in
$D_2(0)\setminus\{0\}$, $w>0$ in $D_2(0)$. Then there exists
$\e_1\in(0,1)$ such that
$$w^{\xi,1-\sqrt\e}(\eta)<u(\eta)$$ for every $\xi,\eta\in\h^n$ and
$\e>0$ such that $\Hn\xi<\e<\e_1$ and $0<\Hn\eta\leq1$.
\end{lem}

\textbf{Proof:} Let $\gamma,\,\e_0\in(0,1)$ be small constants to be
chosen later. Let $\xi,\,\eta\in\h^n$ and $\e>0$ be such that
$\Hn\xi<\e<\e_0$ and $0<\Hn\eta\leq\gamma$. Then one has
$$\Hn{\xi\circ\delta_{1-\sqrt\e}(\eta)}\leq\Hn\xi+(1-\sqrt\e)\Hn\eta<\e+(1-\sqrt\e)\gamma<2$$
and hence in particular $\xi\circ\delta_{1-\sqrt\e}(\eta)\in
D_2(0)$. Moreover for a suitable constant $C>0$ one has
\begin{equation}\label{25}
|\xi\circ\delta_{1-\sqrt\e}(\eta)-\eta|\leq
C\sqrt\e(\sqrt\e+\gamma).
\end{equation}
Since $w\in\mathcal{C}^2\big(D_2(0)\big)$, by relation \eqref{25}
for $\e_0$ and $\delta$ small enough we have that
$$w\big(\xi\circ\delta_{1-\sqrt\e}(\eta)\big)=w(\eta)+O\big(|\xi\circ\delta_{1-\sqrt\e}(\eta)-\eta|\big)=w(\eta)+\sqrt\e
O(\sqrt\e+\gamma),$$ uniformly in $\xi,\,\eta$. Then
\begin{eqnarray}
\nonumber
w^{\xi,1-\sqrt\e}(\eta)-u(\eta)&<&w^{\xi,1-\sqrt\e}(\eta)-w(\eta)\\
\nonumber&=&\left(1-\frac{Q-2}{2}\sqrt\e+O(\e)\right)
w\big(\xi\circ\delta_{1-\sqrt\e}(\eta)\big)-w(\eta)\\
\nonumber&=&\left(1-\frac{Q-2}{2}\sqrt\e+O(\e)\right)
\big(w(\eta)+\sqrt\e O(\sqrt\e+\gamma)\big)-w(\eta)\\
\nonumber&=&\left[-\frac{Q-2}{2}w(\eta)+O(\sqrt\e+\gamma)\right]\sqrt\e,
\end{eqnarray}
uniformly in $\xi,\,\eta$. Thus, for some $\e_0,\,\gamma>0$ small
enough we get
\begin{equation}\label{33}
w^{\xi,1-\sqrt\e}(\eta)<u(\eta)\qquad\text{if
}\,\,\,\,\Hn\xi<\e<\e_0\,\,\,\,\text{ and
}\,\,\,\,0<\Hn\eta\leq\gamma.
\end{equation}
Now let $\xi,\,\eta\in\h^n$ and $\e>0$ be such that $\Hn\xi<\e<\e_0$
and $0<\Hn\eta\leq1$. Then one finds
$$|\xi\circ\delta_{1-\sqrt\e}(\eta)-\eta|\leq
C\sqrt\e.$$ Since $w\in\mathcal{C}^2\big(D_2(0)\big)$ one has
$$w\big(\xi\circ\delta_{1-\sqrt\e}(\eta)\big)=w(\eta)+O\big(|\xi\circ\delta_{1-\sqrt\e}(\eta)-\eta|\big)=w(\eta)+O(\sqrt\e),$$
 uniformly in $\xi,\,\eta$. Choose $\e_1\in(0,\e_0)$ such that
$$O(\sqrt{\e_1})<\min_{\gamma\leq\Hn\zeta\leq1}[u(\zeta)-w(\zeta)].$$
This is possible by our hypotheses on the functions $u,\,w$. Now for
every $\xi,\,\eta\in\h^n$ such that $\Hn\xi<\e<\e_1$ and
$\gamma\leq\Hn\eta\leq1$ we have
\begin{equation}\label{43}
w^{\xi,1-\sqrt\e}(\eta)=w(\eta)+O(\sqrt\e)<w(\eta)+\big(u(\eta)-w(\eta)\big)=u(\eta).
\end{equation}
Relations \eqref{33} and \eqref{43} together give the desired
result. $\qquad\Box$

\begin{lem}\label{2}
Assume $u\in\mathcal{C}^2\big(D_2(0)\setminus\{0\}\big)$,
$w\in\mathcal{C}^2\big(D_2(0)\big)$, $u>w$ in
$D_2(0)\setminus\{0\}$, $w>0$ in $D_2(0)$. Then there exists
$\e_2\in(0,1)$ such that for every $\e\in(0,\e_2)$ one has
$$w^{\xi,\lambda}(\eta)<u(\eta)$$ for
every $\xi\in\h^n$ such that $\Hn\xi<\e$, every $\eta\in\h^n$ such
that $\Hn\eta=1$ and every $\lambda\in[1-\sqrt\e,1+\sqrt\e]$.
\end{lem}

\textbf{Proof:} Let $$\e_0:=\displaystyle\inf_{\partial
D_1(0)}[u(\zeta)-w(\zeta)],$$ then $\e_0>0$. For every
$\e\in(0,1/9)$, every $\lambda\in[1-\sqrt\e,1+\sqrt\e]$ and every
$\xi,\,\eta\in\h^n$ such that $\Hn\xi<\e$ and $\Hn\eta=1$ one has
$\xi\circ\delta_\lambda(\eta)\in D_\frac{3}{2}(0)$. Thus we get
\begin{eqnarray}
\nonumber
w^{\xi,\lambda}(\eta)-u(\eta)\!\!&\!\!=\!\!&\!\!\big[w\big(\xi\circ\delta_\lambda(\eta)\big)-w(\eta)\big]+\big[w(\eta)-u(\eta)\big]
+\big(\lambda^\frac{Q-2}{2}-1\big)\,w\big(\xi\circ\delta_\lambda(\eta)\big)\\
\label{92}\!\!&\!\!\leq\!\!&\!\!\big|w\big(\xi\circ\delta_\lambda(\eta)\big)-w(\eta)\big|-\inf_{\partial
D_1(0)}\big[u(\zeta)-w(\zeta)\big]
+\big|\lambda^\frac{Q-2}{2}-1\big|\sup_{\overline{D_\frac{3}{2}(0)}}w\\
\nonumber\!\!&\!\!\leq\!\!&\!\!\big|w\big(\xi\circ\delta_\lambda(\eta)\big)-w(\eta)\big|-\e_0+C_1\sqrt\e
\end{eqnarray}
for some constant $C_1>0$.

Now notice that for every $\e\in(0,1)$, every
$\lambda\in[1-\sqrt\e,1+\sqrt\e]$ and every $\xi,\,\eta\in\h^n$ such
that $\Hn\xi<\e$ and $\Hn\eta=1$ one has
\begin{equation}\label{44}
|\xi\circ\delta_\lambda(\eta)-\eta|<C_2\sqrt\e
\end{equation}
for some constant $C_2>0$. By the uniform continuity of $w$ on
$\overline{D_\frac{3}{2}(0)}$, since we have
$\eta,\,\big(\xi\circ\delta_\lambda(\eta)\big)\in
\overline{D_\frac{3}{2}(0)}$, by \eqref{44} we can find
$\e_*\in(0,1/9)$ small enough so that for every $\e\in(0,\e_*)$ we
have
$$\big|w\big(\xi\circ\delta_\lambda(\eta)\big)-w(\eta)\big|<\frac{\e_0}{2}.$$
Now let $\e_2=\min\big\{\e_*,\frac{\e^2_0}{16\,C_1^2}\big\}$. By
inequality \eqref{92}, for every $\e\in(0,\e_2)$, every
$\lambda\in[1-\sqrt\e,1+\sqrt\e]$, every $\xi,\,\eta\in\h^n$ such
that $\Hn\xi<\e$ and $\Hn\eta=1$ we have
$$w^{\xi,\lambda}(\eta)-u(\eta)<-\frac{\e_0}{2}+C_1\sqrt\e<-\frac{\e_0}{4}<0.$$
The proof of the lemma is now complete. $\qquad\Box$

\begin{lem}\label{7}
Assume $u\in\mathcal{C}^2\big(D_2(0)\setminus\{0\}\big)$,
$w\in\mathcal{C}^2\big(D_2(0)\big)$, $u>w$ in
$D_2(0)\setminus\{0\}$, $w>0$ in $D_2(0)$. If
$$\liminf_{\xi\rightarrow0}[u(\xi)-w(\xi)]=0,$$
then there exists $\e_3\in(0,1)$ such that
$$\sup_{0<\Hn\eta\leq1}\left\{w^{\xi,1+\frac{\sqrt\e}{2}}(\eta)-u(\eta)\right\}>0$$
for every $\xi\in\h^n$ and every $\e>0$ such that $\Hn\xi<\e<\e_3$.
\end{lem}

\textbf{Proof:} For $\e_3>0$ small enough and every $\e\in(0,\e_3)$,
if $\Hn\xi<\e<\e_3$ one has
\begin{eqnarray}
\nonumber\limsup_{\eta\rightarrow0}\big\{w^{\xi,1+\frac{\sqrt\e}{2}}(\eta)-u(\eta)\big\}&=&w^{\xi,1+\frac{\sqrt\e}{2}}(0)-
\liminf_{\eta\rightarrow0}u(\eta)\\
\label{97}&=&w^{\xi,1+\frac{\sqrt\e}{2}}(0)-w(0)\\
\nonumber&=&\left[1+\frac{Q-2}{2}\,\frac{\sqrt\e}{2}+O(\e)\right]w(\xi)-w(0).
\end{eqnarray}
Now notice that $|\xi|\leq\sqrt2\e$. Since
$w\in\mathcal{C}^2\big(D_2(0)\big)$, we have
$$w(\xi)=w(0)+O(|\xi|)=w(0)+O(\e).$$
Thus by relation \eqref{97} we get
\begin{eqnarray}
\nonumber\limsup_{\eta\rightarrow0}\big\{w^{\xi,1+\frac{\sqrt\e}{2}}(\eta)-u(\eta)\big\}&=&
\left[1+\frac{Q-2}{2}\,\frac{\sqrt\e}{2}+O(\e)\right]\big(w(0)+O(\e)\big)-w(0)\\
\nonumber&=&\left(\frac{Q-2}{4}w(0)+O(\sqrt\e)\right)\sqrt\e>0,
\end{eqnarray}
if $\e\in(0,\e_3)$ and $\e_3$ is small enough. Thus it follows that
$$\sup_{0<\Hn\eta\leq1}\big\{w^{\xi,1+\frac{\sqrt\e}{2}}(\eta)-u(\eta)\big\}>0.\,\qquad\Box$$

We are now ready to prove Theorem \ref{94}.

\textbf{Proof of Theorem \ref{94}:} Without loss of generality, up
to a translation we can assume that $\xi_0=0$.

Arguing by contradiction, suppose that the conclusion of the theorem
does not hold, i.e. that
\begin{equation}\label{99}
\liminf_{\xi\rightarrow0}[u(\xi)-w(\xi)]=0.
\end{equation}
Then let $\e:=\frac{1}{8}\min\{\e_1,\e_2,\e_3\}$, with
$\e_1,\,\e_2,\,\e_3$ being provided by Lemma \ref{15}, Lemma \ref{2}
and Lemma \ref{7} respectively. By Lemma \ref{15} one has
\begin{equation}\label{100}
w^{\xi,1-\sqrt\e}(\eta)<u(\eta)\qquad\text{for all
}\,\eta\in\overline{D_1(0)}\setminus\{0\},
\end{equation}
for every $\xi\in D_\e(0)$. Now for any such $\xi$ define
$$\overline\lambda(\xi):=\sup\left\{\mu\geq1\!-\!\sqrt\e\,\,\Big|\,\,w^{\xi,\lambda}(\eta)<u(\eta)\,\,\,\,\,\forall\,\eta\in\overline{D_1(0)}\setminus\!\{0\},\text{ for all
}\lambda\in[1\!-\!\sqrt\e,\mu]\right\}.$$ By relation \eqref{100} we
have that $\overline\lambda(\xi)$ is well defined and that
$\overline\lambda(\xi)\geq1-\sqrt\e$, for every $\xi\in D_\e(0)$. By
Lemma \ref{7} we also have that
$\overline\lambda(\xi)<1+\frac{\sqrt\e}{2}$ for every $\xi\in
D_\e(0)$.

Now by the definition of $\overline\lambda(\xi)$ one has that
\begin{equation}\label{103}
w^{\xi,\overline\lambda(\xi)}(\eta)\leq u(\eta)\qquad\text{ for
every }\eta\in\overline{D_1(0)}\setminus\!\{0\}\text{ and every
}\xi\in D_\e(0).
\end{equation}
Moreover, since
$\overline\lambda(\xi)\in[1-\sqrt\e,1+\frac{\sqrt\e}{2}]$, by Lemma
\ref{2} one also has that
\begin{equation}\label{104}
w^{\xi,\overline\lambda(\xi)}(\eta)< u(\eta)\qquad\text{ for every
}\eta\in\partial D_1(0)\,\text{ and every }\xi\in D_\e(0).
\end{equation}
We make the following \emph{claim:}
\begin{equation}\label{105}
w^{\xi,\overline\lambda(\xi)}(\eta)< u(\eta)\qquad\text{ for every
}\eta\in\overline{D_1(0)}\setminus\!\{0\}\text{ and every }\xi\in
D_\e(0).
\end{equation}

For every $\xi\in D_\e(0)$,
$\eta\in\overline{D_1(0)}\setminus\!\{0\}$ we have
$\xi\circ\delta_{\overline\lambda(\xi)}(\eta)\in D_\frac{3}{2}(0)$.
By the invariance of the operator $T$, we have that
$$T\big(w^{\xi,\overline\lambda(\xi)},\Hgrad
w^{\xi,\overline\lambda(\xi)},\Hgrad^2w^{\xi,\overline\lambda(\xi)}\big)(\eta)\equiv
T\big(w,\Hgrad
w,\Hgrad^2w\big)\big(\xi\circ\delta_{\overline\lambda(\xi)}(\eta)\big)\leq0\quad\text{
for }\eta\in\overline{D_1(0)}$$ and thus
$$T\big(u,\Hgrad
u,\Hgrad^2u\big)-T\big(w^{\xi,\overline\lambda(\xi)},\Hgrad
w^{\xi,\overline\lambda(\xi)},\Hgrad^2w^{\xi,\overline\lambda(\xi)}\big)\geq0
\qquad\text{ on }\overline{D_1(0)}\setminus\!\{0\}.$$ Now let
$\Omega$ be an open, connected set such that
$\overline\Omega\subset\overline{D_1(0)}\setminus\!\{0\}$. By
condition \eqref{91} we have
$$w^{\xi,\overline\lambda(\xi)}(\eta)-u(\eta)\leq0\qquad\text{ for
every }\eta\in\overline{\Omega}\text{ and every }\xi\in D_\e(0).$$

Recalling the notation $T=T(s,v,U)$ with $s>0$, $v\in\R^{2n}$,
$U\in\mathcal{S}^{2n\times2n}\otimes J\R$, by the regularity of all
the functions involved on $\overline\Omega$, for every $\xi\in
D_\e(0)$ and every $\eta\in\overline\Omega$ we have
\begin{eqnarray}
\nonumber0\!&\!\leq\!&\!T\big(u,\Hgrad
u,\Hgrad^2u\big)-T\big(w^{\xi,\overline\lambda(\xi)},\Hgrad
w^{\xi,\overline\lambda(\xi)},\Hgrad^2w^{\xi,\overline\lambda(\xi)}\big)\\
\nonumber\!&\!=\!&\!\int_0^1\frac{d}{dr}\Big[T\big(ru+(1-r)w^{\xi,\overline\lambda(\xi)},r\Hgrad
u+(1-r)\Hgrad w^{\xi,\overline\lambda(\xi)},\\
\nonumber\!&\!\!&\!\hspace{7.5cm}r\Hgrad^2u+(1-r)\Hgrad^2w^{\xi,\overline\lambda(\xi)}\big)\Big]\,dr\\
\nonumber\!&\!=\!&\!c\big(w^{\xi,\overline\lambda(\xi)}-u\big)+\big\langle
b,\,\Hgrad\big(w^{\xi,\overline\lambda(\xi)}-u\big)\big\rangle_{\R^{2n}}
+\sum_{i,j=1}^{2n}a_{ij}\big[\Hgrad^2\big(w^{\xi,\overline\lambda(\xi)}-u\big)\big]_{ij}\\
\nonumber\!&\!:=\!&\!L\big(w^{\xi,\overline\lambda(\xi)}-u\big),
\end{eqnarray}
with $c\in\mathcal{C}^0(\overline\Omega)$,
$b=(b_1,\ldots,b_{2n})\in\mathcal{C}^0(\overline\Omega,\R^{2n})$,
$a_{ij}\in\mathcal{C}^0(\overline\Omega)$ for every
$i,j=1,\ldots,2n$ being defined by
\begin{eqnarray}
\nonumber c(\eta)\!&\!=\!&\!-\int_0^1\frac{\partial T}{\partial
s}\big(ru(\eta)+(1-r)w^{\xi,\overline\lambda(\xi)}(\eta),r\Hgrad
u(\eta)+(1-r)\Hgrad
w^{\xi,\overline\lambda(\xi)}(\eta),\\
\nonumber \!&\!\!&\!\hspace{7cm}r\Hgrad^2u(\eta)+(1-r)\Hgrad^2w^{\xi,\overline\lambda(\xi)}(\eta)\big)\,dr,\\
\nonumber b_j(\eta)\!&\!=\!&\!-\int_0^1\frac{\partial T}{\partial
v_j}\big(ru(\eta)+(1-r)w^{\xi,\overline\lambda(\xi)}(\eta),r\Hgrad
u(\eta)+(1-r)\Hgrad w^{\xi,\overline\lambda(\xi)}(\eta),\\
\nonumber \!&\!\!&\!\hspace{7cm}r\Hgrad^2u(\eta)+(1-r)\Hgrad^2w^{\xi,\overline\lambda(\xi)}(\eta)\big)\,dr,\\
\nonumber a_{ij}(\eta)\!&\!=\!&\!-\int_0^1\frac{\partial T}{\partial
U_{ij}}\big(ru(\eta)+(1-r)w^{\xi,\overline\lambda(\xi)}(\eta),r\Hgrad
u(\eta)+(1-r)\Hgrad w^{\xi,\overline\lambda(\xi)}(\eta),\\
\nonumber\!&\!\!&\!\hspace{7cm}r\Hgrad^2u(\eta)+(1-r)\Hgrad^2w^{\xi,\overline\lambda(\xi)}(\eta)\big)\,dr.
\end{eqnarray}
By the regularity of $u,\,w$ and since
$\xi\circ\delta_{\overline\lambda(\xi)}(\eta)\in\overline{D_{\frac{3}{2}}(0)}$
for every $\xi\in D_\e(0)$ and every $\eta\in\overline{\Omega}$, we
can find $a,\,b,\,R_1,\,R_2\in\R^+$ with $b>a$ such that
\begin{eqnarray}
\nonumber&&\big(ru+(1-r)w^{\xi,\overline\lambda(\xi)},r\Hgrad
u+(1-r)\Hgrad
w^{\xi,\overline\lambda(\xi)},r\Hgrad^2u+(1-r)\Hgrad^2w^{\xi,\overline\lambda(\xi)}\big)\\
\nonumber&&\hspace{3cm}\in[a,b]\times B_{R_1}(0)\times
B_{R_2}(0)\subset\R^+\times\R^{2n}\times(\mathcal{S}^{2n\times2n}\oplus
J\R)
\end{eqnarray}
for every $r\in[0,1]$, every $\xi\in D_\e(0)$ and every
$\eta\in\overline\Omega$. Thus by our hypotheses on the operator
$T$, see condition \eqref{106}, the matrix
$[a_{ij}(\eta)]_{i,j=1,\ldots,2n}$ is strictly positive definite on
$\overline\Omega$, for every $\xi\in D_\e(0)$.

Since the vector fields $X_1,\ldots,X_n,Y_1,\ldots,Y_n$ satisfy the
H\"{o}rmander condition, see condition \eqref{commutation}, by
Theorem 4.1 in \cite{LanMon} we can conclude that for every $\xi\in
D_1(0)$ we have either
$\big(w^{\xi,\overline\lambda(\xi)}-u\big)\equiv0$ or
$\big(w^{\xi,\overline\lambda(\xi)}-u\big)<0$ in $\overline\Omega$.

Now for any $m\in\N$ we can choose $\Omega=\Omega_m:=D_1(0)\setminus
D_\frac{1}{m}(0)$, and by condition \eqref{104} and the previous
argument we can conclude that
$$\big(w^{\xi,\overline\lambda(\xi)}-u\big)<0\qquad\text{ in }\overline{D_1(0)}\setminus D_\frac{1}{m}(0),\text{ for every }\xi\in
D_\e(0)\text{ and every }m\in\N.$$ Thus our \emph{claim} follows and
\eqref{105} holds.

By the definition of $\overline\lambda(\xi)$ and by relation
\eqref{105}, then one has
\begin{equation}\label{107}
\liminf_{\eta\rightarrow0}[u(\eta)-w^{\xi,\overline\lambda(\xi)}(\eta)]=0\qquad\text{
for every }\xi\in D_\e(0).
\end{equation}
Thus in particular for every $\xi\in D_\e(0)$
\begin{equation}\label{108}
w^{\xi,\overline\lambda(\xi)}(0)=\big(\overline\lambda(\xi)\big)^\frac{Q-2}{2}w(\xi)=\liminf_{\eta\rightarrow0}u(\eta):=\alpha\in(0,\infty).
\end{equation}
By conditions \eqref{105}, \eqref{107} and \eqref{108} and since
$\Delta_{\!_H}u\leq0$, Theorem \ref{93} yields
\begin{equation}\label{109}
V=\Hgrad
w^{\xi,\overline\lambda(\xi)}(0)=\big(\overline\lambda(\xi)\big)^\frac{Q}{2}\Hgrad
w(\xi)\qquad\text{ for every }\xi\in D_\e(0),
\end{equation}
for some fixed $V\in\R^{2n}$. Relations \eqref{108} and \eqref{109}
in turn give $$\Hgrad
w(\xi)=\alpha^{-\frac{Q}{Q-2}}\big(w(\xi)\big)^\frac{Q}{Q-2}V\qquad\text{
for every }\xi\in D_\e(0).$$ Thus for every $\xi\in D_\e(0)$ we get
for $j=1,\ldots,n$
\begin{eqnarray}
\nonumber
X_j^2w(\xi)&=&\frac{Q}{Q-2}\alpha^{-\frac{2Q}{Q-2}}\big(w(\xi)\big)^\frac{Q+2}{Q-2}V_j^2,\\
\nonumber
Y_j^2w(\xi)&=&\frac{Q}{Q-2}\alpha^{-\frac{2Q}{Q-2}}\big(w(\xi)\big)^\frac{Q+2}{Q-2}V_{j+n}^2,
\end{eqnarray}
and hence
$$\Delta_{\!_H}w(\xi)=\frac{Q}{Q-2}\alpha^{-\frac{2Q}{Q-2}}\big(w(\xi)\big)^\frac{Q+2}{Q-2}|V|^2\geq0\qquad\text{ in }D_\e(0).$$
Then it follows that $\Delta_{\!_H}(u-w)\leq0$ in
$D_\e(0)\setminus\{0\}$, with $u-w>0$ in $D_\e(0)\setminus\{0\}$ by
condition (i). Then Lemma \ref{98} yields
$$\inf_{D_\e(0)\setminus\{0\}}(u-w)\geq\inf_{\partial D_\e(0)}(u-w)>0,$$
and thus $$\liminf_{\xi\rightarrow0}[u(\xi)-w(\xi)]>0,$$ which
contradicts our initial assumption, condition \eqref{99}. The proof
of Theorem \ref{94} is now complete. $\qquad\Box$

\section{\textbf{Appendix}}\label{57}

In this section we collect some technical results, which were used
in the course of the proof of Theorem \ref{1} in section \ref{47}.

Moreover we consider the CR map $\check{\varphi}=\varphi\circ\iota$,
which we recall was defined in \eqref{48} by setting
$\check{\varphi}(\xi)=\big(\check{x},\,\check{y},\,\check{t}\big)$
for every $\xi\in\h^n\setminus\{0\}$, with
\begin{equation*}
\check{x}=-\frac{xt+y|z|^2}{\Hn{\xi}^4},\quad\check{y}=\frac{yt-x|z|^2}{\Hn{\xi}^4},\quad\check{t}=\frac{t}{\Hn{\xi}^4},
\end{equation*}
and for every $h=1,\,\ldots,\,n$ we compute
$\Hgrad(\check{x}_h),\,\Hgrad(\check{y}_h),\,\Hgrad(\check{t}),\,\Hgrad^2(\check{x}_h),\,\Hgrad^2(\check{y}_h)$,
$\Hgrad^2(\check{t})$.

\begin{lem}\label{50}
For every $h,j=1,\ldots,n$ we have
\begin{eqnarray}
\nonumber
X_j(\check{x}_h)=-Y_j(\check{y}_h)\!\!&\!\!=\!\!&\!\!-\frac{t}{\Hn{\xi}^4}\,\delta_{\!jh}+\frac{1}{\Hn{\xi}^8}\Big(2\big(|z|^4-t^2\big)(y_h
x_j-x_hy_j)+4t|z|^2(x_jx_h+y_jy_h)\Big)\\
\nonumber
X_j(\check{y}_h)=Y_j(\check{x}_h)\!\!&\!\!=\!\!&\!\!-\frac{|z|^2}{\Hn{\xi}^4}\,\delta_{\!jh}+\frac{1}{\Hn{\xi}^8}\Big(2\big(|z|^4-t^2\big)(x_j
x_h+y_hy_j)+4t|z|^2(y_jx_h-x_jy_h)\Big)\\
\nonumber
X_j(\check{t})\!&\!=\!&\!\frac{1}{\Hn{\xi}^8}\Big(2\big(|z|^4-t^2\big)y_j-4t|z|^2x_j\Big)\\
\nonumber
Y_j(\check{t}\,)\!&\!=\!&\!\frac{1}{\Hn{\xi}^8}\Big(2\big(t^2-|z|^4\big)x_j-4t|z|^2y_j\Big)\\
\nonumber
T(\check{x}_h)\!&\!=\!&\!\frac{1}{\Hn{\xi}^8}\Big(\big(t^2-|z|^4\big)x_h+2t|z|^2y_h\Big)\\
\nonumber
T(\check{y}_h)\!&\!=\!&\!\frac{1}{\Hn{\xi}^8}\Big(\big(|z|^4-t^2\big)y_h+2t|z|^2x_h\Big)\\
\nonumber
T(\check{t})\!&\!=\!&\!\frac{1}{\Hn{\xi}^8}\big(|z|^4-t^2\big)
\end{eqnarray}
where $\delta_{\!jh}$ is the Kronecker's symbol.
\end{lem}
\textbf{Proof:} One has just to use the definition of the vector
fields $T$, $X_j$, $Y_j$ for $j=1,\ldots,n$ given in section
\ref{sec2_1} in order to get the result.$\qquad\Box$

\begin{lem}\label{53} For every $i,h,j=1,\ldots,n$ we have
\begin{eqnarray}
\nonumber
X_jX_i(\check{x}_h)\!&\!=\!&\!\frac{\delta_{ih}}{\Hn{\xi}^8}\big(2(t^2-|z|^4)y_j+4t|z|^2x_j\big)\\
\nonumber\!&\!\!&\!\hspace{0,8cm}-\frac{8}{\Hn{\xi}^{12}}(|z|^2x_j+ty_j)\big(2(|z|^4-t^2)(x_iy_h-y_ix_h)
+4t|z|^2(x_ix_h+y_iy_h)\big)\\
\nonumber\!&\!\!&\!\hspace{0,8cm}+\frac{1}{\Hn{\xi}^8}\Big(8(|z|^2x_j-ty_j)(x_iy_h-y_ix_h)+8(tx_j+|z|^2y_j)(x_ix_h+y_iy_h)\\
\nonumber\!&\!\!&\!\hspace{1,7cm}+\delta_{ij}\big(2(|z|^4-t^2)y_h+4t|z|^2x_h\big)+
\delta_{\!jh}\big(2(t^2-|z|^4)y_i+4t|z|^2x_i\big)\Big)
\end{eqnarray}
\begin{eqnarray}
\nonumber
X_jY_i(\check{x}_h)\!&\!=\!&\!\frac{\delta_{ih}}{\Hn{\xi}^8}\big(2(|z|^4-t^2)x_j+4t|z|^2y_j\big)\\
\nonumber\!&\!\!&\!\hspace{0,8cm}-\frac{8}{\Hn{\xi}^{12}}(|z|^2x_j+ty_j)\big(2(|z|^4-t^2)(y_iy_h+x_ix_h)
+4t|z|^2(y_ix_h-x_iy_h)\big)\\
\nonumber\!&\!\!&\!\hspace{0,8cm}+\frac{1}{\Hn{\xi}^8}\Big(8(|z|^2x_j-ty_j)(y_iy_h+x_ix_h)+8(tx_j+|z|^2y_j)(y_ix_h-x_iy_h)\\
\nonumber\!&\!\!&\!\hspace{1,7cm}+\delta_{ij}\big(2(|z|^4-t^2)x_h-4t|z|^2y_h\big)+
\delta_{\!jh}\big(2(|z|^4-t^2)x_i+4t|z|^2y_i\big)\Big)\\
\nonumber
Y_jX_i(\check{x}_h)\!&\!=\!&\!\frac{\delta_{ih}}{\Hn{\xi}^8}\big(2(|z|^4-t^2)x_j+4t|z|^2y_j\big)\\
\nonumber\!&\!\!&\!\hspace{0,8cm}+\frac{8}{\Hn{\xi}^{12}}(|z|^2y_j-tx_j)\big(2(|z|^4-t^2)(x_hy_i-x_iy_h)
-4t|z|^2(y_iy_h+x_ix_h)\big)\\
\nonumber\!&\!\!&\!\hspace{0,8cm}-\frac{1}{\Hn{\xi}^8}\Big(8(|z|^2y_j+tx_j)(y_ix_h-x_iy_h)+8(|z|^2x_j-ty_j)(y_iy_h+x_ix_h)\\
\nonumber\!&\!\!&\!\hspace{1,7cm}+\delta_{ij}\big(2(|z|^4-t^2)x_h-4t|z|^2y_h\big)+
\delta_{\!jh}\big(2(t^2-|z|^4)x_i-4t|z|^2y_i\big)\Big)\\
\nonumber
Y_jY_i(\check{x}_h)\!&\!=\!&\!\frac{\delta_{ih}}{\Hn{\xi}^8}\big(2(|z|^4-t^2)y_j-4t|z|^2x_j\big)\\
\nonumber\!&\!\!&\!\hspace{0,8cm}-\frac{8}{\Hn{\xi}^{12}}(|z|^2y_j-tx_j)\big(2(|z|^4-t^2)(y_iy_h+x_ix_h)
+4t|z|^2(y_ix_h-x_iy_h)\big)\\
\nonumber\!&\!\!&\!\hspace{0,8cm}+\frac{1}{\Hn{\xi}^8}\Big(8(|z|^2y_j+tx_j)(y_iy_h+x_ix_h)+8(|z|^2x_j-ty_j)(x_iy_h-y_ix_h)\\
\nonumber\!&\!\!&\!\hspace{1,7cm}+\delta_{ij}\big(2(|z|^4-t^2)y_h+4t|z|^2x_h\big)+
\delta_{\!jh}\big(2(|z|^4-t^2)y_i-4t|z|^2x_i\big)\Big)
\end{eqnarray}
Moreover
\begin{eqnarray}
\nonumber X_jX_i(\check{y}_h)=X_jY_i(\check{x}_h)&\qquad&Y_jY_i(\check{y}_h)=-Y_jX_i(\check{x}_h),\\
\nonumber
Y_jX_i(\check{y}_h)=Y_jY_i(\check{x}_h)&\qquad&X_jY_i(\check{y}_h)=-X_jX_i(\check{x}_h).
\end{eqnarray}

Finally
\begin{eqnarray}
\nonumber X_jX_i(\check{t})\!&\!=\!&\!-\frac{8}{\Hn{\xi}^{12}}(|z|^2x_j+ty_j)\big(2(|z|^4-t^2)y_i-4t|z|^2x_i\big)\\
\nonumber\!&\!\!&\!\hspace{0,8cm}+\frac{1}{\Hn{\xi}^8}\Big(8(|z|^2x_j-ty_j)y_i-8(tx_j+|z|^2y_j)x_i
-4t|z|^2\delta_{ij}\Big)\\
\nonumber
X_jY_i(\check{t})\!&\!=\!&\!-\frac{8}{\Hn{\xi}^{12}}(|z|^2x_j+ty_j)\big(2(t^2-|z|^4)x_i-4t|z|^2y_i\big)\\
\nonumber\!&\!\!&\!\hspace{0,8cm}+\frac{1}{\Hn{\xi}^8}\Big(8(ty_j-|z|^2x_j)x_i-8(tx_j+|z|^2y_j)y_i
+2(t^2-|z|^4)\delta_{ij}\Big)\\
\nonumber
Y_jY_i(\check{t})\!&\!=\!&\!-\frac{8}{\Hn{\xi}^{12}}(|z|^2y_j-tx_j)\big(2(t^2-|z|^4)x_i-4t|z|^2y_i\big)\\
\nonumber\!&\!\!&\!\hspace{0,8cm}+\frac{1}{\Hn{\xi}^8}\Big(8(|z|^2x_j-ty_j)y_i-8(|z|^2y_j+tx_j)x_i
-4t|z|^2\delta_{ij}\Big).
\end{eqnarray}
\end{lem}
\textbf{Proof:} One needs only use the definition of the vector
fields $T$, $X_j$, $Y_j$ for $j=1,\ldots,n$ and lemma \ref{50} to
conclude.$\qquad\Box$

\begin{rem}
By Lemma \ref{50} and Lemma \ref{53}, for every $h=1,\ldots,n$ we
have
$$\Hgrad(\check{x}_h)=-J\Hgrad(\check{y}_h),\qquad\Hgrad^2(\check{x}_h)=-J\Hgrad^2(\check{y}_h)$$
in $\h^n\setminus\{0\}$, where $J\in\text{Mat}(2n,\R)$ is the matrix
defined in \eqref{45}. Moreover
\begin{eqnarray}
\nonumber|\Hgrad(\check{x}_h)|&=&\left(\sum_{j=1}^n\Big(|X_j(\check{x}_h)|^2+|Y_j(\check{x}_h)|^2\Big)\right)^\frac{1}{2}\,\,=\,\,\frac{1}{\Hn{\xi}^2},\\
\nonumber|\Hgrad(\check{y}_h)|&=&\left(\sum_{j=1}^n\Big(|X_j(\check{y}_h)|^2+|Y_j(\check{y}_h)|^2\Big)\right)^\frac{1}{2}\,\,=\,\,\frac{1}{\Hn{\xi}^2}
\end{eqnarray}
for every $\xi\in\h^n\setminus\{0\}$.
\end{rem}

\begin{lem}\label{52} Let $v=(p,q)\in\R^n\times\R^n$, $v\neq0$, $s>0$ and let
$\xi_0:=(x_0,y_0,t_0)$ with
$$x_0=-\frac{(Q-2)s}{|v|^2}p,\quad y_0=-\frac{(Q-2)s}{|v|^2}q,\quad t_0=0.$$
Define $z_0=(x_0,y_0)$ and
$$\lambda:=\Hn{\xi_0}=\big(|x_0|^2+|y_0|^2\big)^\frac{1}{2}=\frac{(Q-2)s}{|v|}.$$
Let
\begin{equation*}
\psi(\xi)=\varphi\circ\iota\circ\delta_{\lambda^{-2}}(\xi)=
\big(\lambda^2\check{x},\,\lambda^2\check{y},\,\lambda^4\check{t}\big)
\end{equation*}
for every every $\xi=(x,y,t)\in\h^n\setminus\{0\}$, and consider
also a positive function $\phi\in\mathcal{C}^\infty(\h^n)$ such that
$\phi(\xi_0)=s$, $\Hgrad\phi(\xi_0)=v$ and $\Hgrad^2\phi(\xi_0)=U$.
Then one has $\phi_\psi(\psi^{-1}\big(\xi_0)\big)=s$,
$\Hgrad\phi_\psi(\psi^{-1}\big(\xi_0)\big)=0$ and
\begin{eqnarray}
\nonumber\Hgrad^2\phi_\psi(\psi^{-1}(\xi_0))\!\!&\!\!=\!\!&\!\!G\bigg[-\frac{Q}{Q-2}s^{-1}Jv\otimes
    Jv+\frac{2}{Q-2}s^{-1}v\otimes v+\frac{1}{Q-2}s^{-1}|v|^2I_{2n}\\
\nonumber\!&\!\!&\!\hspace{0,5cm}+J^TUJ+\frac{4}{|v|^4}\Big(\crochet{v}{Uv}_{\R^{2n}}Jv\otimes
    Jv+\crochet{Jv}{UJv}_{\R^{2n}}v\otimes v\\
\label{54}\!&\!\!&\!\hspace{0,5cm}-\crochet{Jv}{Uv}_{\R^{2n}}v\otimes
    Jv-\crochet{v}{UJv}_{\R^{2n}}Jv\otimes v\Big)+\frac{2}{|v|^2}\Big(Jv\otimes
    J^TU^Tv\\
\nonumber\!&\!\!&\!\hspace{0,5cm}-v\otimes J^TU^TJv+J^TUv\otimes
   Jv-J^TUJv\otimes v\Big)\bigg]G,
\end{eqnarray}
with $G$, $J$ being defined as in \eqref{45}.
\end{lem}
\textbf{Proof:} Notice that
$\psi=\check{\varphi}\circ\delta_{\lambda^{-2}}$, so that by
formulae \eqref{4} and \eqref{48} and by Proposition \ref{51}, for
every $\xi\in\h^n\setminus\{0\}$ one has
\begin{eqnarray}
\nonumber\big(\phi_{\check{\varphi}}\big)_{\delta_{\lambda^{-2}}}(\xi)&=&
\lambda^{-(Q-2)}\phi_{\check{\varphi}}(\lambda^{-2}x,\lambda^{-2}y,\lambda^{-4}t)\\
\label{71}&=&\left(\frac{\lambda}{\Hn{\xi}}\right)^{Q-2}\phi(\lambda^2\check{x},\lambda^2\check{y},\lambda^4\check{t})\\
\nonumber&=&\phi_\psi(\xi).
\end{eqnarray}
Now notice also that for every $\xi\in\h^n\setminus\{0\}$ one has
$\psi^2(\xi)=\xi$, so that $\psi^{-1}(\xi)=\psi(\xi)$ for all
$\xi\in\h^n\setminus\{0\}$. Then by the definition of $\lambda$ and
$\xi_0$ in particular one has
$$\psi^{-1}(\xi_0)\,\,=\,\,\psi(\xi_0)\,\,=\,\,(-y_0,-x_0,0).$$ Thus using formula
\eqref{71} one immediately gets
$$\phi_\psi\big(\psi^{-1}(\xi_0)\big)\,\,=\,\,\phi(\xi_0)\,\,=\,\,s$$

By formulae \eqref{71} and \eqref{grad_tr}, for every
$\xi=(x,y,t)\in\h^n\setminus\{0\}$ one also has
\begin{eqnarray}
\label{65}&&\Hgrad{\phi_\psi}(\xi)=\Hgrad\Big(\big(\phi_{\check{\varphi}}\big)_{\delta_{\lambda^{-2}}}\Big)(\xi)=
\lambda^{-Q}\Hgrad\phi_{\check{\varphi}}\big(\lambda^{-2}x,\lambda^{-2}y,\lambda^{-4}t\big).
\end{eqnarray}
Evaluating equality \eqref{65} at $\psi^{-1}(\xi_0)$ and using
formula \eqref{62} we get
\begin{eqnarray}
\nonumber&&\Hgrad{\phi_\psi}\big(\psi^{-1}(\xi_0)\big)=
\frac{(Q-2)\lambda^{Q-2}}{|z_0|^Q}u(\xi_0)\left(
                                            \begin{array}{c}
                                              y_0 \\
                                              x_0 \\
                                            \end{array}
                                          \right)
+\frac{\lambda^Q}{|z_0|^Q}\,E\cdot\Hgrad{u}(\xi_0),
\end{eqnarray}
where the matrix $E$, which was defined in \eqref{63}, is now
evaluated at the point $\psi^{-1}(\xi_0)$. Recalling the definition
of $x_0$ and $y_0$, from the previous equality one concludes that
$\Hgrad{\phi_\psi}\big(\psi^{-1}(\xi_0)\big)=0$ as claimed.

In a similar way, by formulae \eqref{Hess_tr} one also has
\begin{eqnarray}
\label{74}&&\Hgrad^2\phi_\psi(\xi)=\Hgrad^2\Big(\big(\phi_{\check{\varphi}}\big)_{\delta_{\lambda^{-2}}}\Big)(\xi)=
\lambda^{-(Q+2)}\Hgrad^2\phi_{\check{\varphi}}\big(\lambda^{-2}x,\lambda^{-2}y,\lambda^{-4}t\big),
\end{eqnarray}
so that
\begin{eqnarray}
\nonumber&&\Hgrad^2\phi_\psi(\psi^{-1}(\xi_0))=
\lambda^{-(Q+2)}\Hgrad^2\phi_{\check{\varphi}}\big(-\lambda^{-2}y_0,-\lambda^{-2}x_0,0\big).
\end{eqnarray}
Now in order to get formula \eqref{54}, and thus conclude the proof
of the lemma, it's sufficient to substitute in the previous equality
the expression for $\Hgrad^2\phi_{\check{\varphi}}$ which is
provided by formula \eqref{67}, where all the functions appearing
there are to be evaluated at the point
$\xi=(-\lambda^{-2}y_0,-\lambda^{-2}x_0,0)$. $\qquad\Box$

\begin{lem}\label{72} Let $\lambda>0$ and let
\begin{eqnarray}
\nonumber&&\xi_0\,\,=\,\,(x_0,y_0,t_0)\,\,=\,\,(0,0,\lambda^2),\\
\nonumber&&\psi(\xi)\,\,=\,\,\varphi\circ\iota\circ\delta_{\lambda^{-2}}(\xi)\,\,=\,\,
\big(\lambda^2\check{x},\,\lambda^2\check{y},\,\lambda^4\check{t}\big)
\end{eqnarray}
for every every $\xi=(x,y,t)\in\h^n\setminus\{0\}$. Consider a
positive function $\phi\in\mathcal{C}^\infty(\h^n)$ such that
$\phi(\xi_0)=1$, $\Hgrad\phi(\xi_0)=0$ and $\Hgrad^2\phi(\xi_0)=U$.
Then one has $\phi_\psi\big(\psi^{-1}(\xi_0)\big)=1$,
$\Hgrad\phi_\psi\big(\psi^{-1}(\xi_0)\big)=0$ and
\begin{eqnarray}
\nonumber\Hgrad^2\phi_\psi\big(\psi^{-1}(\xi_0)\big)\!\!&\!\!=\!\!&\!\!-\frac{Q-2}{\lambda^2}J+GUG,
\end{eqnarray}
with $G$, $J$ being defined as in \eqref{45}.
\end{lem}
\textbf{Proof:} One has $\psi(\xi_0)=\psi^{-1}(\xi_0)=\xi_0$ and
$\Hn{\xi_0}=\lambda$. Thus, using formula \eqref{71} one immediately
gets $$\phi_\psi\big(\psi^{-1}(\xi_0)\big)=\phi(\xi_0)=1.$$ Next,
evaluating equality \eqref{65} at $\xi_0=\psi^{-1}(\xi_0)$ and using
formula \eqref{62} yields
$$\Hgrad\phi_\psi\big(\psi^{-1}(\xi_0)\big)=0.$$
Finally, we evaluate \eqref{74} at the point
$\xi_0=\psi^{-1}(\xi_0)$ and we use formula \eqref{67}. Since in
this case we have $E=-G$, where $E$ is the matrix defined in
\eqref{63} now evaluated at the point $(0,0,\lambda^{-2})$, we get
$$\Hgrad^2\phi_\psi\big(\psi^{-1}(\xi_0)\big)=-\frac{Q-2}{\lambda^2}J+GUG,$$
and the proof of the lemma is now complete. $\qquad\Box$

One can prove the following lemma in the same way.

\begin{lem}\label{73} Let $\lambda>0$ and let
\begin{eqnarray}
\nonumber&&\xi_0\,\,=\,\,(x_0,y_0,t_0)\,\,=\,\,(0,0,-\lambda^2),\\
\nonumber&&\psi(\xi)\,\,=\,\,\varphi\circ\iota\circ\delta_{\lambda^{-2}}(\xi)\,\,=\,\,
\big(\lambda^2\check{x},\,\lambda^2\check{y},\,\lambda^4\check{t}\big)
\end{eqnarray}
for every every $\xi=(x,y,t)\in\h^n\setminus\{0\}$. Consider a
positive function $\phi\in\mathcal{C}^\infty(\h^n)$ such that
$\phi(\xi_0)=1$, $\Hgrad\phi(\xi_0)=0$ and $\Hgrad^2\phi(\xi_0)=U$.
Then one has $\phi_\psi\big(\psi^{-1}(\xi_0)\big)=1$,
$\Hgrad\phi_\psi\big(\psi^{-1}(\xi_0)\big)=0$ and
\begin{eqnarray}
\nonumber\Hgrad^2\phi_\psi\big(\psi^{-1}(\xi_0)\big)\!\!&\!\!=\!\!&\!\!\frac{Q-2}{\lambda^2}J+GUG,
\end{eqnarray}
with $G$, $J$ being defined as in \eqref{45}.
\end{lem}

The following lemma is a consequence of Bony's Strong Maximum
Principle, see \cite{Bony}.

\begin{lem}\label{98}
Let
$u\in\mathcal{C}^2\big(\overline{D_r(\xi_0)}\setminus\{\xi_0\}\big)$
for some $r>0$ and $\xi_0\in\h^n$. Assume
\begin{itemize}
\item[i)] $\Delta_{\!_H}u\leq0\,\text{ on
}\,D_r(\xi_0)\setminus\{\xi_0\}$,
\item[ii)] $u^-(\xi)=o\left(\frac{1}{\Hn{\xi-\xi_0}^{Q-2}}\right)$ as $\,\xi\rightarrow\xi_0$, where
$u^-:=-\min\{u,0\}$.
\end{itemize}
Then
$\displaystyle\inf_{D_r(\xi_0)\setminus\{\xi_0\}}u\geq\inf_{\partial
D_r(\xi_0)} u$.
\end{lem}

\textbf{Proof:} Up to a translation $\tau_{\xi_0}$ and a dilation
$\delta_r$, we can assume without loss of generality that
$D_r(\xi_0)\equiv D_1(0)$.

Now consider the function $w(\xi)=\frac{1}{\Hn{\xi}^{Q-2}}$ on
$\h^n\setminus\{0\}$, which is a multiple of the fundamental
solution of $-\Delta_{\!_H}$ on $\h^n$, centered at $0\in\h^n$. Let
$\displaystyle m_0:=\min_{\partial D_1(0)}u$ and define
$v:=\frac{m_0-u}{w}$ on $\overline{D_1(0)}\setminus\{0\}$. Then for
every $\e>0$ one has
\begin{equation*}
\Delta_{\!_H}v+\frac{2}{w}\crochet{\Hgrad w}{\Hgrad
v}_{\R^{2n}}=-\frac{\Delta_{\!_H}u}{w}\geq0
\end{equation*}
on $D_1(0)\setminus D_\e(0)$. Since $w\in
C^\infty(\h^n\setminus\{0\})$ and $w>0$ in $\h^n\setminus\{0\}$, by
the Strong Maximum Principle proved by Bony in \cite{Bony} one has
that
\begin{equation}\label{101}
\sup_{\overline{D_1(0)\setminus D_\e(0)}} v\leq\sup_{\partial
D_1(0)\cup\partial D_\e(0)}v^+=\sup_{\partial D_\e(0)}v^+,
\end{equation}
since by definition one has $v\leq0$ on $\partial D_1(0)$. For every
fixed $\xi\in D_1(0)\setminus\{0\}$, one has that $\xi\in
D_1(0)\setminus D_\e(0)$ for every $\e>0$ small enough. Then using
\eqref{101} one obtains
\begin{equation}\label{102}
v(\xi)\leq\sup_{\partial D_\e(0)} v^+\leq\sup_{\partial
D_\e(0)}\frac{c_0+u^-}{w}\leq c_0\e^{Q-2}+\sup_{\partial
D_\e(0)}\frac{u^-}{w}
\end{equation}
for every $\e>0$ small enough. By condition $ii)$, one has that
$\frac{u^-(\eta)}{w(\eta)}$ tends to $0$ as $\eta\rightarrow0$, and
hence $\displaystyle\sup_{\partial D_\e(0)}\frac{u^-}{w}$ also tends
to $0$ as $\e\rightarrow0$. Then passing to the limit as
$\e\rightarrow0$ in \eqref{102} yields $v(\xi)\leq0$, for every
$\xi\in D_1(0)\setminus\{0\}$. By the definition of $v$ thus we have
$$u(\xi)\geq c_0=\min_{\partial D_1(0)}u\qquad\text{for every
}\xi\in\overline{D_1(0)}\setminus\{0\}.\qquad\Box$$

\vspace{0,2cm}
\begin{center}
{\bf Acknowledgements}
\end{center}

The authors would like to thank Vittorio Martino for his useful
comments, which helped in the exposition of the results.

The second author wishes to thank the Department of Mathematics of
Rutgers University for the hospitality and stimulating environment
provided during the preparation of this work. The second author also
wishes to thank the Nonlinear Analysis Center at Rutgers for having
hosted his visit.

\bibliographystyle{plain}

\bibliography{bib_Li_Monticelli}

\end{document}